\documentclass[final,5p,times,twocolumn]{elsarticle}

 \usepackage{graphicx}
 \usepackage{color}

\usepackage{amssymb}
\usepackage{amsmath}
 \usepackage{amsthm}

\usepackage{array}

\biboptions{sort&compress}

\journal{Engineering Analysis with Boundary Elements (Accepted June 11, 2019)}

\newcommand{\Bmath}[1]{\boldsymbol{#1}}
\newcommand{\Vector}[1]{\Bmath{#1}}
\newcommand{\Vr}{\Vector{r}}
\newcommand{\x}{\Vector{r}}
\newcommand{\xx}{\Vector{r}_i}
\newcommand{\nablaxx}{\nabla\!_i}
\newcommand{\Vn}{\Vector{n}}
\newcommand{\Ver}{\Vector{e}_r}

\newcommand{\Int}[4]{\int_{#1}^{#2}\! #3\,d{#4}\,}
\newcommand{\OInt}[4]{\oint_{#1}^{#2} #3\,d{#4}\,}
\newcommand{\Gint}[2]{\Int{\Gamma_{#1}}{}{#2}{\Gamma}}
\newcommand{\GOint}[2]{\OInt{\Gamma_{#1}}{}{#2}{\Gamma}}

\newcommand{\NodeA}{\alpha}
\newcommand{\NodeB}{\beta}
\newcommand{\NodeC}{\gamma}
\newcommand{\NodeCC}{\rm c}
\newcommand{\NodeAdomI}{{\NodeA(1)}}
\newcommand{\NodeAdomII}{{\NodeA(2)}}

\newcommand{\NodeBdomI}{{\NodeB(1)}}
\newcommand{\NodeBdomII}{{\NodeB(2)}}

\newcommand{\NodeCdomI}{{\NodeC(1)}}
\newcommand{\NodeCdomII}{{\NodeC(2)}}
\newcommand{\EleA}{{\rm A}}
\newcommand{\EleB}{{\rm B}}

\newcommand{\OMint}[2]{\Int{\Omega_{#1}}{}{\,\,#2}{\Omega}}
\newcommand{\ROUND}[2]{\frac{\partial{#1}}{\partial{#2}}}
\newcommand{\RoundR}[1]{\ROUND{#1}{r}}
\newcommand{\IP}{\hspace*{-0.0833em}\cdot\hspace*{-0.0833em}}
\newcommand{\Grad}[1]{\nabla #1}
\newcommand{\Gradxx}[1]{\nablaxx #1}
\newcommand{\Tan}{\tau}

\newcommand{\Rint}[3]{\Int{#1}{#2}{#3}{r}}
\newcommand{\TypeAB}{\gamma}
\newcommand{\Typeal}{\TypeAB}

\newcommand{\sal}{\Tan_\Typeal}
\newcommand{\nal}{n_\Typeal}
\newcommand{\Vsal}{\Vector{\Tan}_\Typeal}
\newcommand{\Vnal}{\Vector{n}_\Typeal}
\newcommand{\Dyad}[1]{\overleftrightarrow{{\bf{#1}}}}
\newcommand{\DCi}{\Dyad{C_{\mathit i}}}
\newcommand{\DCiinv}{\Dyad{C_{\mathit i}^{\rm{-1}}}}

\newcommand{\Diff}[1]{\underset{\TypeAB\,:\,\EleA-\EleB}{{\rm Diff}}\left[\mathstrut{#1}\right]}
\newcommand{\FuncAt}[2]{\left.{#1}\right|_{#2}}

\newcommand{\Dtheta}{\varDelta\theta}

\newcommand{\NormAlt}[1]{%
  \settoheight{\dimen0}{$\mathstrut #1$}%
  \advance\dimen0 .0pt%
  \settodepth{\dimen1}{$#1$}%
  \advance\dimen1 4pt%
  \,\hbox{
    \vrule height\dimen0 depth\dimen1\,%
    \vrule height\dimen0 depth\dimen1\,%
    \hbox{$#1$}%
    \,%
    \vrule height\dimen0 depth\dimen1\,%
    \vrule height\dimen0 depth\dimen1\,%
  }\,%
} 
\newcommand{\Norm}[1]{\NormAlt{#1}}

\newcommand{\Define}{\triangleq}

\newcommand{\EQ}[1]{Eq.~(\ref{#1})}
\newcommand{\EQtwo}[2]{Eqs.~(\ref{#1}) and (\ref{#2})}
\newcommand{\EQthree}[3]{Eqs.~(\ref{#1}), (\ref{#2}) and (\ref{#3})}

\newcommand{\Fig}[1]{Fig.\!\,\,\ref{fig:#1}}

\begin{document}

\begin{frontmatter}

\title{%
Desingularization of matrix equations employing hypersingular integrals in
boundary element methods using double nodes
}%

\cortext[cor1]{Corresponding author}
\author{Satoshi Tomioka\corref{cor1}}
\ead{tom@qe.eng.hokudai.ac.jp}
\author{Shusuke Nishiyama\corref{}}
\author{Yutaka Matsumoto\corref{}}
\author{Naoki Miyamoto\corref{}}

\address{Faculty of Engineering, Hokkaido University, Sapporo, 060-8628, Japan}

\begin{abstract}
%%%%%%%%%%%%%%%%%%%%%%%%%%%%%%%%%%%%%%%
In boundary element methods,
using double nodes at corners is a useful approach to uniquely
define the normal direction of boundary elements.
However, matrix equations constructed by conventional boundary integral equations (CBIEs)
become singular
under certain combinations of double node boundary conditions.
In this paper, we analyze the singular conditions of the CBIE formulation
for cases where the boundary conditions at the double node are imposed
by combinations of Dirichlet, Neumann, Robin, and interface conditions. 
To address this singularity,
we propose the use of hypersingular integral equations (HBIEs)
for wave propagation problems that obey the Helmholtz equation.
To demonstrate the applicability of HBIE, we compare three types of simultaneous equations:
(i) CBIE, (ii) partial-HBIE where the
HBIE is only applied to the double nodes at corners while the CBIE is applied to the other nodes,
and (iii) full-HBIE where the HBIE is applied to all nodes.
Based on our numerical results, we observe the following results.
The singularity of the matrix equations
for problems with any combination of boundary conditions
can be resolved by both full-HBIEs and partial-HBIEs, and the
partial-HBIE exhibits better accuracy than the full-HBIE.
Furthermore, the computational cost of partial-HBIEs is smaller than that of full-HBIEs.

\end{abstract}


\end{frontmatter}

%%%%%%%%%%%%%%%%%%%%%%%%%%%%%%%%%%%%%%%%%%%%%%%%%%%%%%%%%%%%%%%%%%

\section{Introduction}
\label{sec:intro}
The boundary element method (BEM), the finite difference method (FDM), and the finite element method (FEM)
have been commonly used to solve boundary value problems.
In the BEM, a set of simultaneous equations for determining unknown variables at nodes on the boundary
is constructed in discretized boundary integral equations.
The variables in simultaneous equations 
are nodal field values and normal derivatives only on individual boundary nodes,
whereas the variables in the FDM or FEM are field values at domain nodes which are placed in the entire domain enclosed by the boundary.
Therefore, the number of variables in the BEM is much smaller than that in the FDM or FEM,
which is one of the advantages of BEM.
Furthermore, the BEM can be easily applied to external problems, such as wave scattering problems,
since it does not require the placement of nodes in a domain spreading to infinity
and it does not require any other boundary conditions to represent radiation at a boundary enclosing the domain considered.

In most boundary problems,
the field values, $u$, along the boundary are continuous;
however, the normal derivatives, $q$, are discontinuous at corners
since the normal directions at any corner point are different.
By using a linear element or higher-order elements,
the boundary elements share the nodes at both ends of the element with adjacent boundary elements.
In this case, the normal direction, $\Vn$, at the corner node
cannot be defined uniquely
since the single node at any corner belongs to two boundary elements
with different normal directions.

There are two approaches to addressing the definition problem of the normal direction.
The first approach involves the use of non-conforming elements, which are also called discontinuous elements.
In the non-conforming element,
collocation nodes that represent $u$ and $q$ do not coincide with geometric nodes,
but they do so in the conforming element.
The non-conforming element has been applied to several problems;
e.g.,
elastostatic problems \cite{Manolis:1986,Parreira:1988,Olukoko:1993,Blazquez:1994,Huesmann:1994,Paris:1995,Blazquez:1998},
fluid flow problems \cite{Patterson:1982,Dyka:1989},
and acoustic problems \cite{Silva:1993}.
Although the accuracy between the non-conforming element
and the conforming element were compared
by Manolis and Banerjee \cite{Manolis:1986}, and Parreira \cite{Parreira:1988},
they arrived at different conclusions.

The second approach includes a double node technique \cite{Brebbia:1984:Sec_5_2}
or a multiple node technique for three-dimensional problems.
In the double node technique,
two normal derivatives, $q_\NodeA=\Vn_\NodeA\IP\Grad{u}$ and $q_\NodeB=\Vn_\NodeB\IP\Grad{u}$, and a field, $u$,
are defined at the corner node,
where $\Vn_\NodeA$ and $\Vn_\NodeB$ denote
the directions normal to the two boundary elements connected to the corner node.
However,
a set of simultaneous equations, called a matrix equation, becomes singular
under certain boundary conditions;
i.e., the rank of the matrix equation is reduced
since some node equations are redundant. Further details will be presented in Sec.~\ref{sec:CBIE-Rank_deficient}.

To address the rank reduction problem caused by the double nodes,
there are two categories of approaches.
The first category involves the use of a local relation for each double node \cite{Walker:1989,Yan:1994,Kassab:1994,Gao:2000}.
By using a Taylor expansion around the corner node,
this relation is described as a linear combination of the two normal derivatives and the field values at neighbor nodes of the corner node.
The local relation is replaced by one of the redundant equations that reduces the rank;
therefore, the matrix equation still includes a square matrix.
In the second category,
extra node equations are employed \cite{Mitra:1987,Mitra:1993,Subia:1995,Deng:2013,Zheng:2018}.
Mitra and Ingber \cite{Mitra:1987} proposed a technique
for replacing one of the redundant equations in each corner by an extra node equation
with respect to an extra collocation node placed outside the domain considered.
Following this study, the authors mentioned that ``external collocation yields a coefficient matrix with a large number'' \cite{Mitra:1993}, and
they improved the method using the extra node equations
so that the location of the extra node on the boundary elements connects to the double node \cite{Mitra:1993,Subia:1995}.
Subia, Ingber, and Mitra demonstrated that
there are no significant differences in accuracy
between the method of the extra node equation and the non-conforming method \cite{Subia:1995}.
The method that uses the extra node equations was extended to the problems of interface boundaries
at which two or more domains are connected \cite{Deng:2013,Zheng:2018}.
Using these methods, the number of variables is not increased
since the field at the extra node is known.
Therefore, the coefficient matrix is the square matrix, which is similar to the matrix of the first category.
If we simply add the local relations
shown in the first category or the extra node equations in the second category
instead of replacing them,
the number of equations becomes larger than the number of variables,
which is referred to as an overdetermined problem.
We can solve this equation by using least-square methods,
but the computational time required to solve this equation is much greater than solving the general simultaneous equation;
therefore, the replacements are generally applied.
The replacement of the equation should be performed individually
while examining the types of boundary conditions at the corner.
The individual examination increases the complexity of the programming of widely applicable BEM codes that include many types of boundary conditions.

In addition to the corner node problems, there are rank deficient or large condition number problems.
We focus on two problems related to hypersingular boundary integral equations.
The first is a non-uniqueness or a spurious solution problem 
in an external field for a wave scattering problem that obeys the Helmholtz equation,
in which
the domain considered is outside of the boundary enclosing a scatterer.
In this situation, spurious solutions can be obtained at the eigenfrequency of
the scatterer.
To resolve this problem,
two major approaches were proposed.
Schenck employed equations related to additional nodes
in the scatterer region and the method is called
`Combined Helmholtz Integral Equation Formulation (CHIEF)' \cite{Schenck:1968}.
Chen et~al.~\cite{IL_Chen:2001} applied the additional node equations to interior problems
in which additional points are placed outside of the boundary.
These methods are similar to the extra node equation methods for the corner problems
shown in the previous paragraph; however, they require a solver based on the least-square method
since the set of simultaneous equations becomes an overdetermined equation.
The other approach is referred to as the Burton-Miller method \cite{Burton-Miller:1971,Benthien:1997,Diwan:2013,Langrenne:2015}.
Burton and Miller \cite{Burton-Miller:1971} represented a boundary integral equation (BIE) for each node
using a linear combination of two types of BIE called a conventional BIE (CBIE) and a hypersingular BIE (HBIE).
In the CBIE,
the field $u$ at the field point is denoted by integrals over the boundary
on which sources are distributed.
The HBIE is obtained by taking a gradient of the CBIE.
The singularity of the integral in the HBIE is stronger than that of the CBIE;
therefore, it is called a hypersingular integral.
Bentihien and Schenck \cite{Benthien:1997} reviewed the non-uniqueness problem with comparisons
of other methods, including the CHIEF and Burton-Miller methods.

The other rank deficient problem is found at a degenerate boundary which appears either at a crack in an elastostatic problem \cite{Portela:1992,JT_Chen:1994,Chyuan:2003,Lu:2010}
or at both surfaces of a thin metal with zero-thickness in an electromagnetic problem \cite{Chyuan:2003}.
To resolve these problems, the CBIEs are applied to one side of the crack or the thin metal,
and the HBIEs are applied to the other side.
These methods are referred to as dual-BEM.

As shown in the previous two paragraphs, the use of HBIEs is effective for resolving rank deficient problems.
In this paper, we will demonstrate that
the corner singular problem for wave propagation problems can be solved
by only using node equations based on HBIEs.
In addition, to suppress the rank deficiency caused by the corners,
we do not need to use HBIEs for every node.
We will also illustrate that only the replacement of the node equations related to the corner node by HBIEs is sufficient,
which is similar to the aforementioned dual-BEM.

In the application of HBIEs, the regularization of hypersingularity is a key issue.
The authors developed 
an analytical regularization of the two-dimensional Helmholtz equation \cite{Tomioka:2010};
the other regularization methods are also found in the references of that study.
In the regularization of HBIEs, we include the relations between two normal derivatives at the double nodes,
similar to the local relation methods described above.
Therefore, the method by HBIE can be considered an extension that uses the local relation.
However,
the contributions to the double node from all the nodes are considered in the method by HBIE;
whereas the local relation method represents the relations between local nodal quantities only.

The method proposed in this paper to overcome the corner problem caused by the double nodes does not
require local relations at corners,
extra node equations,
or least-square methods.
Our method can also be applied to any kind of boundary conditions of boundary elements
that include corner nodes.
In addition, the proposed method can also be applied to interface boundary conditions.
From these characteristics,
no additional effort is required in prepossessing to prepare the input data for solving the boundary value problem.

The outline of this paper is as follows.
In Sec.~\ref{sec:CBIE-Rank_deficient},
we demonstrate why the rank of the coefficient matrix of the CBIE
is reduced in the case where the double node is employed. We also demonstrate
the condition that results in rank deficiency based on the nature of the discretized node equations.
In Sec.~\ref{sec:HBIE}, we illustrate why the HBIE does not cause a rank deficiency.
The numerical results and discussions for simple waveguide problems are presented in Sec.~\ref{sec:results}
to demonstrate that the HBIE is applicable to any combination of boundary conditions.
The advantages of the partial-HBIE method
in which the HBIE is applied to
only the double nodes and the CBIE is applied to other nodes are also shown.
Finally, the summary is presented in Sec.~\ref{sec:conclusion}.

\section{Rank deficiency problem in CBIEs}
\label{sec:CBIE-Rank_deficient}
To illustrate rank deficient conditions for
a set of discretized node equations of CBIEs,
the discretization using linear elements is first shown;
the double node technique that defines two sub-nodes for a double node are then shown;
and lastly,
the rank deficiency conditions are discussed
based on comparisons between two equations for the two sub-nodes.

\subsection{Discretization of CBIEs}
\label{sec:CBIE-discretization}
A complex-valued time harmonic scalar wave
satisfies the following Helmholtz equation
with an assumed time factor ${\rm e}^{{\rm j}\omega t}$,
where $\rm j$ is an imaginary unit and $\omega$ is the angular frequency:
\begin{align}
  \label{Helmholtz}
  &&\nabla^2 u(\x)+k^2 u(\x)=0, \qquad \x\in\Omega,
\end{align}
where $k$ indicates the wave number, which is a ratio of $\omega$ to the wave velocity,
and $\Omega$ represents the spatial domain
considered.
A fundamental solution $u^*(\x,\xx)$,
which 
represents the contribution to a field point $\xx$ from a unit source placed
at a source point $\x$ in free space
satisfies
\begin{align}
  \label{funda1}
  \nabla^2 u^*(\x,\xx)+k^2 u^*(\x,\xx)=-\delta(\x-\xx),
\end{align}
where the differential operator $\nabla$ operates only on $\x$,
but not on $\xx$.
This equation can be solved analytically.
In two-dimensional problems,
the outward propagating wave that obeys \EQ{funda1} is 
\begin{align}
   u^*(\x,\xx)=\frac{1}{4\rm j}H_0^{(2)}(kr),\qquad r=|\x-\xx|,
\end{align}
where the function $H_0^{(2)}(kr)$ is a second kind 0-th order Hankel function.

Using Green's second identity and some integral operations for \EQtwo{Helmholtz}{funda1},
we obtain the conventional boundary integral equation (CBIE),
\begin{align}
  \label{bie1}
    c(\xx)&\,u(\xx)
     =\GOint{}{\left[ u^*(\x,\xx)\Grad{u(\x)}\IP\Vn
       %     \right.\nonumber\\&&\left.\
             -u(\x) \Grad{ u^*(\x,\xx)}\IP\Vn \right]},
          \nonumber\\
    %=\GOint{}{\left[ u^*\,\Grad{u}\IP\Vn
    %         -u\,\Grad{ u^*}\IP\Vn \right]},
\end{align}
where $\Gamma$ denotes a boundary surrounding $\Omega$,
$\x$ is the position of the source point on the boundary,
$\xx$ is the position of the field point,
$\Vn$ is the outward-pointing normal unit vector;
the contour integral is evaluated as a Cauchy principal value,
and $c(\xx)$ is the result of the following evaluation
of Dirac's delta function:
\begin{align}
  \label{def-c}
  c(\xx)u(\xx)
  \Define\OMint{}{u(\x)\delta(\x-\xx)}
  =\OMint{}{\delta(\x-\xx)}\,u(\xx)\,.
\end{align}
The coefficient $c(\xx)$ depends on
the shape of the boundary $\Gamma$ at the field point $\xx$.
Because of the nature of Dirac's delta function,
when $\xx$ is located inside and outside the domain, $c(\xx)$ evaluates to 1 and 0, respectively.
In the case where $\xx$ is located on $\Gamma$,
$c(\xx)$ is equal to the ratio of the interior angle $\Dtheta$ to the whole angle;
e.g., $\Dtheta/2\pi$ for 2-dimensional problems.
In addition,
when we introduce a fundamental solution to a Laplace equation, $u_L^*$,
which satisfies
\begin{align}
  \label{funda-laplace}
  \nabla^2 u_L^*(\x,\xx)=-\delta(\x-\xx),
\end{align}
the coefficient $c(\xx)$ can be expressed by a boundary integral:
\begin{align}
  \label{equi-pot}
  c(\xx)=\OMint{}{\delta(\x-\xx)}
   =-\GOint{}{\Grad u_L^*\IP\Vn},
\end{align}
which is called an equipotential condition \cite{Brebbia:1980}.

When $\xx$ is located on $\Gamma$ and the boundary is partitioned into a number of boundary elements,
\EQ{bie1} can be written as a node equation for the node $i$ as a discrete algebraic expression:
\begin{align}
  \label{CBIE-point_i}
  &
  \sum_{j\in I}\left(c_i\delta_i^j+h_i^j\right)u_j^{}-\sum_{j\in I} g_i^j q_j^{}=0,
\end{align}
where both the superscripts and the subscripts $i$ and $j$
denote the node numbers and not the element numbers;
$I$ denotes a set of the node numbers, which includes $N$ members;
$\delta_i^j$ denotes the Kronecker's delta;
and $q\equiv\Grad u\IP\Vn$.
The factors $g_i^j$ and $h_i^j$ are results of boundary integrals
in \EQ{bie1} as shown below.

In \EQ{CBIE-point_i}, there are two variables, $u_j$ and $q_j$, at $\x_j$.
One of them is specified by a boundary condition as a boundary value which is denoted by $\overline{x}_j$, and the other is an unknown variable denoted $x_j$; therefore, the number of variables is equal to that of the nodes, $N$.
Since the node $i$ can be placed at every boundary node,
we can obtain $N$ equations of \EQ{CBIE-point_i} as
\begin{align}
  \label{System_eq-point_i}
  \sum_{j\in I} a_i^j x_j=\overline{b}_i, \quad \overline{b}_i=\sum_{j\in I} b_i^j \overline{x}_j
  \qquad \mbox{ for each }i\in I,
\end{align}
where
both $a_i^j$ and $b_i^j$ are either $h_i^j$ or $g_i^j$ when simple boundary conditions are specified.
The matrix composed of $a_i^j$ is called a coefficient matrix in the following discussions.

The coefficient $g_i^j$ represents the contribution to $u_i$
from the nodal quantity $q_j$ at the node $j$,
and $h_i^j$ represents the contribution to $u_i$ from $u_j$.
These coefficients are evaluated by boundary integrals,
which depend on a method that interpolates $u(\Vr)$ and $q(\Vr)$ along the boundary element.
To examine the nature of $g_i^j$ and $h_i^j$,
we present evaluations using a shape function to interpolate them.
Let us consider the discretization of the boundary integral of $u^*(\x,\xx)\,q(\x)$.
By denoting the $k$-th boundary element as $\Gamma_k$ and a local node number as $l$,
$q(\x)$ on $\Gamma_k$ is represented
by a linear combination of shape functions $\phi^{(k,l)}(\x)$ and
$q$ at the boundary nodes:
\begin{align}
  \label{linear-interporation}
  q(\x)=\sum_{l=1}^{N_l} \phi^{(k,l)}(\x) q_{(k,l)}^{}\quad \mbox{ on\ }\Gamma_k,
\end{align}
where $N_l$ denotes the number of nodes in a boundary element;
e.g., $N_l=2$ for the linear element, $N_l=3$ for the second-order element.
A global node number $j$ can be mapped from the local node numbers $(k,l)$
by a permutation matrix $m_{(k,l)}^j$:
\begin{align}
   q_{(k,l)}^{}=\sum_j m^j_{(k,l)} q_j,
\end{align}
where $m^j_{(k,l)}$ has the value 1 if $(k,l)$ and $j$ are associated, and 0 otherwise.
Although the notation of $m^j_{(k,l)}$ is used here, it provides a symbolic meaning, and
it is treated as a mapping function such as $j=m(k,l)$ in actual coding
to avoid summation procedures with respect to $k$.
By using these definitions,
the boundary integral is discretized as
\begin{align}
  \label{use-permutation}
  &
  \GOint{}{u^*(\x,\xx)q(\x)}
   \nonumber\\&\quad
   =\sum_k
    \Gint{k}{u^*(\x,\xx)\left(\sum_l\phi^{(k,l)}(\x) \left(\sum_j m^j_{(k,l)}q_j\right)\right)}
   \nonumber\\&\quad
   %\nonumber\\&\qquad
   =\sum_j\left(\sum_k\sum_l m_{(k,l)}^j g^{(k,l)}_{i}\right) q_j
   =\sum_j g_i^j q_j,
\end{align}
where $g_i^{(k,l)}$ and $g_i^j$ are defined as
\begin{align}
  &
   \label{g_int}
   g^{(k,l)}_{i}=\Gint{k}{u^*(\x,\xx)\phi^{(k,l)}(\x)},
  \\
  &
   \label{g_int2}
   g_i^j=\sum_k\sum_l  m^j_{(k,l)}  g^{(k,l)}_{i}.
\end{align}
The integral on the right-hand side in \EQ{g_int} is evaluated by a numerical integral
such as Gauss quadrature or by an analytical integral in the case of $\xx\in\Gamma_k$,
which is called a singular integral.
Similarly, the other integral in \EQ{bie1} is evaluated as
\begin{align}
  &
  \GOint{}{q^*(\x,\xx)u(\x)}
   =\sum_j h_i^j u_j,
   \\
  &
  \label{h_int}
   h_i^{(k,l)}=\Gint{k}{q^*(\x,\xx)\phi^{(k,l)}(\x)},
  \\
  &
   \label{h_int2}
   h_i^j=\sum_k\sum_l m^j_{(k,l)} h_i^{(k,l)},
\end{align}
where $q^*(\x,\xx)=\Grad u^*(\x,\xx)\IP\Vn(\x)$, which is not the derivative at $\xx$; i.e., $q^*(\x,\xx)$ does not depend on $\Vn(\xx)$.

To discuss the rank of the matrix equation shown in \EQ{System_eq-point_i}, the dependencies of both $g_i^j$ and $h_i^j$
on $\Vn(\xx)$ which is the normal unit vector at $\xx$ are important.
First, $g_i^j$ is independent of $\Vn(\xx)$
since \EQtwo{g_int}{g_int2} do not include $\Vn(\xx)$
regardless of $j$ or $k$.
In contrast,
the case of $h_i^j$ is different from that of $g_i^j$.
Although $h_i^j$ in \EQtwo{h_int}{h_int2} looks independent of $\Vn(\xx)$
even when $q^*(\x,\xx)$ is expressed using $\Vn(\x)$,
there is an exception at $i=j$ in which $\Vn(\x)$ becomes $\Vn(\xx)$.
The dependency of this exceptional case is presented in the next section.

%%%%%%%%%%%%%%%%%%%%%%%%%%%%%%%%%%%%%%%%%%%%%%%%%%%%%%%%

\begin{figure}[tb]
  \centering
  \includegraphics[width=0.6\hsize]{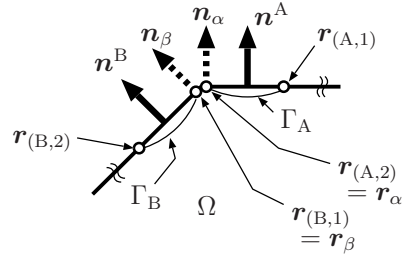}
  \caption{%
    Configuration of double node.
    Thick lines depict boundary elements.
    The two circles which are in contact express sub-nodes which represent a double node at a corner.
    The distance between the sub-nodes are zero; i.e., they are connected to a zero-sized element.
    Each sub-node belonging to the double node is connected to one non-zero-sized element.
    The outward-pointing normal vector is defined for each boundary element
    as $\Vn^\EleA$ or $\Vn^\EleB$.
    The normal vector at each sub-node, $\Vn_\NodeA$ or $\Vn_\NodeB$, can be defined uniquely
    as $\Vn_\NodeA=\Vn^\EleA$ and $\Vn_\NodeB=\Vn^\EleB$.
  }
  \label{fig:double_node}
\end{figure}

\subsection{Double node technique}
\label{sec:double_node}
In the node equation for node $i$ shown in \EQ{CBIE-point_i},
both $u_j$ and $q_j$ are values of the nodes,
which are located at the two ends of the boundary element in the case of the linear element.
When the node $\Vr_i$ is located at a corner,
the normal derivative $q_i$ of the linear element cannot be defined uniquely since 
the node belongs to two elements with different normal directions.
Using a double node technique is one of the solutions.

Figure \ref{fig:double_node} illustrates a configuration of a double node.
The node $\Vr_i$ is represented by two sub-nodes at the same position;
$\Vr_{\NodeA}=\Vr_{\NodeB}$.
The node at $\Vr_{\NodeA}$ is connected to both $\Gamma_{\EleA}$ and a zero-sized element
between $\Vr_{\NodeA}$ and $\Vr_{\NodeB}$, and vice versa.
The normal direction of the sub-nodes $\NodeA$ and $\NodeB$
can be determined from the directions normal to $\Gamma_{\EleA}$ and $\Gamma_{\EleB}$, respectively;
therefore, the direction of the derivatives of $q_{\NodeA}$ and $q_{\NodeB}$ can be defined individually.
Since the integral along the zero-sized element is identically zero regardless of the normal direction,
we do not need to define the normal direction for the zero-sized element.
Although the node $\NodeA$ does not belong to $\Gamma_\EleB$,
contributions of $g^{(\EleB,l)}_{\NodeA}$ and $h^{(\EleB,l)}_{\NodeA}$,
which are evaluated by analytical integrals as singular integrals for accurate evaluation,
are similar to $g^{(\EleA,l)}_{\NodeA}$ and $h^{(\EleA,l)}_{\NodeA}$.

Both $u$ and $q$ are defined at each sub-node in the same way as
ordinary nodes, and the boundary condition is imposed for each sub-node.
Either $u_\NodeA$ or $u_\NodeB$ is unnecessary since $u_\NodeA=u_\NodeB$;
however, the same representation as the ordinary nodes simplifies the programming effort.
In this case,
the number of unknown variables increases
by the number of double nodes, which equals the number of corners. 
Since one double node $i$ is replaced by the two sub-nodes $\NodeA$ and $\NodeB$ for each corner,
the number of node equations is also increased by the number of corners.
Consequently,
a set of CBIEs can be expressed by a matrix equation with a square coefficient matrix
even when we apply double nodes.

Since the integral along the zero-sized elements is zero,
the right-hand sides of \EQ{equi-pot} for $i=\NodeA$ and $i=\NodeB$
are the same.
Therefore, the coefficients $c_{\NodeA}$ and $c_{\NodeB}$ must have the same value: 
\begin{align}
  \label{c-ia-ib}
  c_{\NodeA}^{}=c_{\NodeB}^{}.
\end{align}

In the case of linear elements,
the singular integral of $h_i^j$ for the corners ($i,j=\NodeA,\NodeB$) becomes zero
since the $\x-\xx$ for $i=\NodeA$ or $\NodeB$ is perpendicular to $\Vn_j$ for $j=\NodeA$ or $\NodeB$:
\begin{align}
  \label{h-ia-ib-sing}
  & h_{\NodeA}^{\NodeA}=h_{\NodeB}^{\NodeA}
   =h_{\NodeA}^{\NodeB}=h_{\NodeB}^{\NodeB}
   =0&& \mbox{ for linear elements}.
\end{align}
Therefore, summarizing the discussions
in the last paragraph in Sec.~\ref{sec:CBIE-discretization} and this result,
we obtain the relations for the sub-nodes as
\begin{align}
  \label{g-ia-ib}
  & g_{\NodeA}^j=g_{\NodeB}^j&& \mbox{ for any $j$},
   \\
  \label{h-ia-ib-nosing}
  & h_{\NodeA}^j=h_{\NodeB}^j&& \mbox{ for any $j$}.
\end{align}

\subsection{Rank deficiency conditions in CBIEs}
\label{sec:Rank_deficient_condition-CBIE}
By using the double node technique,
the boundary nodes including sub-nodes
are defined at the ends of boundary elements, and 
the boundary conditions are defined at each node.
In problems with a single medium,
there are three common types of boundary conditions;
Dirichlet condition, Neumann condition, and Robin condition.
These conditions are given, respectively, as
\begin{align}
  \label{Dirichlet}
  &u_j=\overline{u}_j &&\mbox{ for } j\in I_D,
  \\
  \label{Neumann}
  &q_j=\overline{q}_j &&\mbox{ for } j\in I_N,
  \\
  \label{Robin}
  &\kappa_u u_j+\kappa_q q_j=\overline{\psi}_j&&\mbox{ for }j\in I_R,
\end{align}
where the over-bars denote values imposed by the boundary conditions,
$\kappa_u$ and $\kappa_q$ are the given constants determined by the problem considered,
and the sets of the node numbers with corresponding boundary conditions are
denoted by $I_D$, $I_N$, and $I_R$, respectively.

When we consider a multi-media problem,
there exists a condition at the interface between the media.
The interface condition is expressed by two continuous conditions for
$u$ and $q$.
The continuous conditions are given as
\begin{align}
  \label{Interface-cond}
  &u_{j^{(2)}}=u_{j^{(1)}},
   \quad
   q_{j^{(2)}}=-\kappa_{21} q_{j^{(1)}}
  &&\mbox{ for }j^{(1)},j^{(2)}\in I_I,
\end{align}
where media numbers are denoted by (1) and (2);
the nodes $j^{(1)}$ and $j^{(2)}$ are located at the same positions;
$\kappa_{21}$ is determined by media constants of $\Omega^{(1)}$ and $\Omega^{(2)}$;
and $I_I$ denotes the set of nodes with interface conditions.

In this section,
we first present three simple examples when two sub-nodes of a double node at a corner point
belong to $I_D$ or $I_N$,
and discuss why the rank deficient problem arises.
Then, the case in which the sub-nodes 
belong to $I_R$ or $I_I$ are presented.

The combinations of the boundary conditions
where the two sub-nodes, $\NodeA$ and $\NodeB$, belong to $I_D$ or $I_N$
are classified into the following three cases:
\begin{itemize}
  \item both include the Dirichlet conditions:
    %\\\quad
       $\NodeA\in I_D$ and $\NodeB\in I_D$,
  \item both include the Neumann conditions:
    %\\\quad
      $\NodeA\in I_N$ and $\NodeB\in I_N$,
  \item one includes the Dirichlet condition and the other includes the Neumann condition:
    \\\quad
      $\NodeA\in I_D$ and $\NodeB\in I_N$, or $\NodeA\in I_N$ and $\NodeB\in I_D$.
\end{itemize}

\subsubsection{Case of two Dirichlet conditions}
\label{sec:Dirichlet-case-CBIE}
Under these conditions ($\NodeA\in I_D$ and $\NodeB\in I_D$),
the sub-node equations of \EQ{CBIE-point_i} for $i={\NodeA}$ and ${\NodeB}$
can be arranged so that the terms including known values are moved to the right-hand side and
the terms related to $\NodeA$ and $\NodeB$ are moved out from the summation as
\begin{align}
  \label{Both-Dirichlet-ia}
  &
  -g_{\NodeA}^{\NodeA}q_{\NodeA}^{}-g_{\NodeA}^{\NodeB}q_{\NodeB}^{}
  % +\sum_{j\in {I{\setminus\{\EleA,\EleB\}}}}
  +\hspace*{-0.5em}\sum_{j\in I\setminus\{\NodeA,\NodeB\}}\hspace*{-0.5em}
    a_{\NodeA}^j x_{j}
  =\overline{b}_{\NodeA},
  \\
  \label{Both-Dirichlet-ib}
  &
  -g_{\NodeB}^{\NodeA}q_{\NodeA}^{}-g_{\NodeB}^{\NodeB}q_{\NodeB}^{}
  % +\sum_{j\in{I{\setminus\{\EleA,\EleB\}}}}
  +\hspace*{-0.5em}\sum_{j\in I\setminus\{\NodeA,\NodeB\}}\hspace*{-0.5em}
    a_{\NodeB}^j x_{j}
  =\overline{b}_{\NodeB},
\end{align}
where $I=I_D\cup I_N$, i.e., all nodes,
and $I{\setminus\{\NodeA,\NodeB\}}$ denotes a set of all node numbers except $\NodeA$ and $\NodeB$;
and 
the detailed descriptions of the third terms on the left-hand sides and the right-hand sides are written as
\begin{align}
  \label{Both-Dirichlet-ax}
   &
   \hspace*{-0.5em}\sum_{j\in I\setminus\{\NodeA,\NodeB\}}\hspace*{-0.5em}
    a_{i}^j x_{j}
    =\sum_{j\in I_N}h_{i}^j u_j
    -\hspace*{-0.75em}\sum_{{j\in I_D},{j\ne{\NodeA,\NodeB}}}\hspace*{-0.75em} g_{i}^j q_j,
  \\
  \label{Both-Dirichlet-b}
    &\overline{b}_{i}
      =
      -\sum_{j\in I_D}\left(c_{i}^{}\delta_{i}^j+h_{i}^j\right)\overline{u}_j
      +\sum_{j\in I_N}g_{i}^j \overline{q}_j.
\end{align}
Comparing the coefficients of the terms on the left-hand side of the two sub-node equations
shown in \EQtwo{Both-Dirichlet-ia}{Both-Dirichlet-ib},
we can evaluate whether the rank of the coefficient matrix is reduced or not.
From \EQtwo{g-ia-ib}{h-ia-ib-nosing},
the contributions $h_i^j$ and $g_i^j$ ($i=\NodeA,\NodeB$) to the node $j$ are the same.
The coefficients $c_\alpha$ and $c_\beta$ share the same value from \EQ{c-ia-ib};
however, the multiplied terms, $u_j$,
are different since $u_j$ is also multiplied by the Kronecker's delta in \EQ{CBIE-point_i}.
Therefore,
we can examine rank reduction by examining which terms simply include the Kronecker's delta.
The Kronecker's delta is not found
in the first and the second terms
on the left-hand sides in \EQtwo{Both-Dirichlet-ia}{Both-Dirichlet-ib}.
For the third terms, node $j$ does not include $\NodeA$ and $\NodeB$;
therefore, the $\delta_\NodeA^j$ and $\delta_\NodeB^j$ are not included.
They only appear in the right-hand sides, which are not related to rank reduction.
Therefore, the left-hand sides of the two equations are identical,
and the rank is always reduced by these two sub-node equations;
i.e., the coefficient matrix becomes a singular matrix.
 
\subsubsection{Case of two Neumann conditions}
\label{sec:Neumann-case-CBIE}
As in the previous sub-section,
two sub-node equations of \EQ{CBIE-point_i} in the case of $\NodeA\in I_N$ and $\NodeB\in I_N$
can be arranged as
\begin{align}
  \label{Both-Neumann-ia}
  \left(c_{\NodeA}^{}+h_{\NodeA}^{\NodeA}\right)u_{\NodeA}^{}
  +h_{\NodeA}^{\NodeB}u_{\NodeB}^{}
    &
  %+\sum_{j\in {I{\setminus\{\NodeA,\NodeB\}}}}
  +\hspace*{-0.5em}\sum_{j\in I\setminus\{\NodeA,\NodeB\}}\hspace*{-0.5em}
    a_{\NodeA}^j x_{j}
  =\overline{b}_{\NodeA},
  \\
  \label{Both-Neumann-ib}
  h_{\NodeB}^{\NodeA}u_{\NodeA}^{}
  +\left(c_{\NodeB}^{}+h_{\NodeB}^{\NodeB}\right)u_{\NodeB}^{}
    &
  %+\sum_{j\in {I{\setminus\{\NodeA,\NodeB\}}}}
  +\hspace*{-0.5em}\sum_{j\in I\setminus\{\NodeA,\NodeB\}}\hspace*{-0.5em}
    a_{\NodeB}^j x_{j}
  =\overline{b}_{\NodeB}.
\end{align}
Although the definitions of $a_i^j$ and $\overline{b}_i$ are different
from \EQtwo{Both-Dirichlet-ax}{Both-Dirichlet-b},
the third terms of \EQtwo{Both-Neumann-ia}{Both-Neumann-ib}
are identical since they do not include the Kronecker's delta.
In contrast, the first and second terms are different.
By applying \EQtwo{c-ia-ib}{h-ia-ib-sing},
the coefficients associated with $u_\NodeA$ and $u_\NodeB$ in \EQ{Both-Neumann-ia} are $c_\NodeA$ and 0, respectively;
while those in \EQ{Both-Neumann-ib} are 0 and $c_\NodeA$, respectively.
Therefore, the rank of the matrix that includes the sub-node equations is not reduced
in the case of $\NodeA\in I_N$ and $\NodeB\in I_N$.
In addition, the right-hand sides are the same since they do not include the Kronecker's delta.
This means that the results of these two sub-node equations involve the relation $u_\NodeA=u_\NodeB$.

\subsubsection{Case of coupled Dirichlet and Neumann conditions}
\label{sec:Coupled-case-CBIE}
In the case of $\NodeA\in I_D$ and $\NodeB\in I_N$,
two sub-node equations are
\begin{align}
  \label{Dirichlet-Neumann-ia}
  -g_{\NodeA}^{\NodeA}q_{\NodeA}^{}
  &
  +h_{\NodeA}^{\NodeB}u_{\NodeB}^{}
  %+\sum_{j\in {I{\setminus\{\NodeA,\NodeB\}}}}
  +\hspace*{-0.5em}\sum_{j\in I\setminus\{\NodeA,\NodeB\}}\hspace*{-0.5em}
    a_{\NodeA}^j x_{j}
  =\overline{b}_{\NodeA},
  \\
  \label{Dirichlet-Neumann-ib}
  -g_{\NodeB}^{\NodeA}q_{\NodeA}^{}
  &
  +\left(c_{\NodeB}^{}+h_{\NodeB}^{\NodeB}\right)u_{\NodeB}^{}
  +\hspace*{-0.5em}\sum_{j\in I\setminus\{\NodeA,\NodeB\}}\hspace*{-0.5em}
    a_{\NodeB}^j x_{j}
  =\overline{b}_{\NodeB}.
\end{align}
Since the second terms on the left-hand sides of the above equations are different,
the rank of the coefficient matrix is not reduced
in the case of $\NodeA\in I_D$ and $\NodeB\in I_N$.

According to the examples in Secs.~\ref{sec:Neumann-case-CBIE} and ~\ref{sec:Coupled-case-CBIE},
we can understand that
the two sub-node equations for $\NodeA$ and $\NodeB$ are different
when either $c_{\NodeA}$ or $c_{\NodeB}$ is included in
the coefficients of unknown variables in the two sub-node equations.

\subsubsection{Case of Robin condition}
\label{Robin-case-CBIE}
We consider the case in which the Robin conditions are imposed on at least one of the two sub-nodes.
Eliminating $u_j$ from \EQ{CBIE-point_i} using \EQ{Robin},
and arranging the equation,
we obtain:
\begin{align}
  \label{disc-Robin}
  &
  \sum_{j\in I_N}\left(c_i^{}\delta_i^j+h_i^j\right)u_j
  -\sum_{j\in I_D} g_i^j q_j
  %\nonumber\\&
  %\qquad
  -\sum_{j\in I_R}\left(\frac{\kappa_q}{\kappa_u}\left(c_i^{}\delta^j_i+h_i^j\right)+g_i^j\right)q_j
  =\overline{b}_i.
\end{align}
When $i=\NodeA\in I_R$, the term with the factor $c_{\NodeA}$ for unknown $q_{\NodeA}$ remains
in the third terms of the left-hand side.
This satisfies the condition
described at the end of Sec.~\ref{sec:Coupled-case-CBIE},
which ensures that the equations of the two sub-nodes
are independent and not identical.
It does not depend on the type of boundary condition of the node $\NodeB$.

%%%
\begin{figure}[tb]
  \centering
  \includegraphics[width=0.9\hsize]{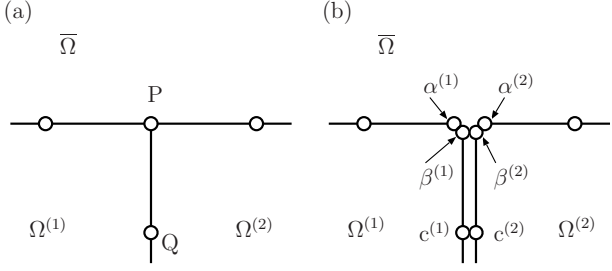}
  \caption{%
    Duplicated node belonging to boundary elements with interface condition.
    (a) Original positions of nodes.
    (b) Definition of sub-nodes.
    Two domains $\Omega^{(1)}$ and $\Omega^{(2)}$ are 
    in contact at an interface.
    The domain $\overline{\Omega}$ is outside of the domains considered.
    The node {\rm Q} in (a) is a general node
    belonging to the boundary element with the interface condition
    that can be separated by two nodes $\NodeCC^{(1)}$ and $\NodeCC^{(2)}$ in (b)
    where the positions are the same.
    Similar to {\rm Q}, the node {\rm P} is a node belonging to the interface;
    however, it is also
    the double node.
    Consequently, {\rm P} is split into four sub-nodes;
    $\NodeA^{(1)}$, $\NodeB^{(1)}\in\Omega_1$, and $\NodeA^{(2)}$, $\NodeB^{(2)}\in\Omega_2$.
    For each pair of
    $\left(\NodeB^{(1)},\NodeB^{(2)}\right)$ or $\left(\NodeCC^{(1)},\NodeCC^{(2)}\right)$,
    the interface conditions are satisfied, and both $u$ and $q$ at these nodes are unknown.
    Either a Dirichlet, Neumann, or Robin condition is imposed
    at $\NodeA^{(1)}$ and $\NodeA^{(2)}$.
  }
  \label{fig:double_node-interface}
\end{figure}
%%%
\begin{figure}
  \centering
  \includegraphics[width=0.9\hsize]{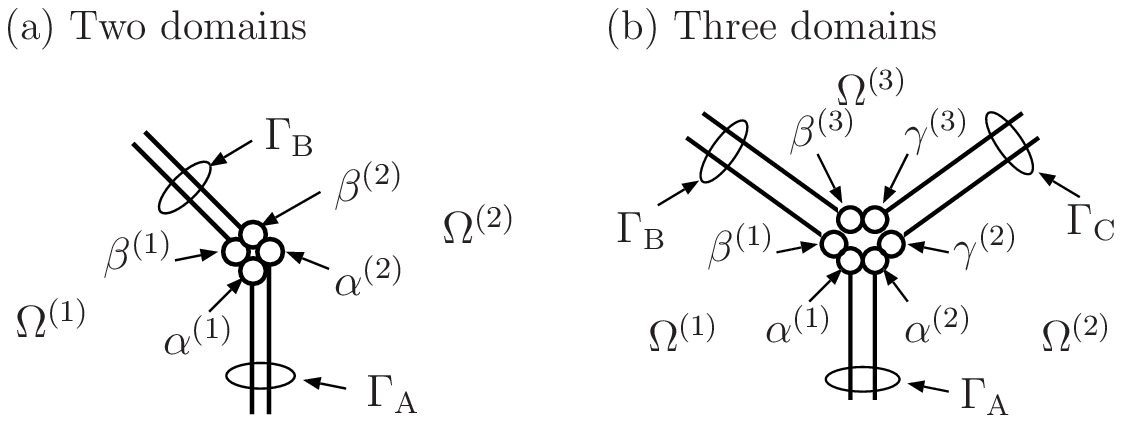}
  \caption{%
    Definition of sub-nodes at an intersection of the interface boundary elements.
    Between the pairs of sub-nodes,
    $(\NodeAdomI,\NodeAdomII)$,
    $(\NodeBdomI,\NodeBdomII)$,
    and $(\NodeCdomI,\NodeCdomII)$,
    the interface boundary conditions are imposed.
  }
  \label{fig:double_node-interface-domain}
\end{figure}

\subsubsection{Case of interface condition}
\label{sec:interface-case-CBIE}
Figure \ref{fig:double_node-interface} illustrates a configuration around a double node
where two domains are in contact.
The corner node labeled {\rm P} is represented by the four sub-nodes at the corner;
two sub-nodes labeled $\NodeA^{(1)}$ and $\NodeA^{(2)}$ do not belong to the interface boundary,
and the other two sub-nodes labeled $\NodeB^{(1)}$ and $\NodeB^{(2)}$ belong to the interface boundary.
There are eight quantities related to these four sub-nodes.
Of the eight quantities, two quantities at the sub-nodes $\NodeA^{(1)}$ and $\NodeA^{(2)}$
are considered ordinary boundary conditions,
and two variables at $\NodeB^{(1)}$ and $\NodeB^{(2)}$ can be eliminated by the interface conditions
shown in \EQ{Interface-cond}, which, for $\NodeB^{(1)}$ and $\NodeB^{(2)}$, are
rewritten as
\begin{align}
  \label{u-cont}
  &u_{\NodeBdomII}^{}=u_{\NodeBdomI}^{},
  %\\
  \quad
  %\label{q-cont}
  q_{\NodeBdomII}^{}=-\kappa_{21} q_{\NodeBdomI}^{},
\end{align}
where the double suffixes, such as $u_{\NodeB^{(1)}}$, are denoted
by $u_{\NodeBdomI}$ for simplicity.
Consequently,
there are four unknown variables and four equations with respect to the double node
when
we apply these conditions before constructing the simultaneous equations.

As mentioned in the previous sub-sections,
the coefficients unrelated to the sub-nodes do not affect rank reduction;
therefore, we only consider the nature of the sub-matrix
composed of the coefficients for the four unknown variables.

In the case of $\NodeA^{(1)}\in I_D$ and $\NodeA^{(2)}\in I_D$,
the given values are $\overline{u}_{\NodeAdomI}$ and $\overline{u}_{\NodeAdomII}$,
and the unknown variables are
$q_{\NodeAdomI}$, $u_{\NodeBdomI}$, $q_{\NodeBdomI}$, $q_{\NodeAdomII}$, $u_{\NodeBdomII}$, and $q_{\NodeBdomII}$.
The independent variables are reduced by using \EQ{u-cont}
such as
$\{q_{\NodeAdomI}, q_{\NodeBdomI},q_{\NodeAdomII},u_{\NodeBdomI}\}$.
The terms associated with these independent variables in the individual sub-node equations 
for $i=\NodeA^{(1)}$, $\NodeB^{(1)}$, $\NodeA^{(2)}$, and $\NodeB^{(2)}$, 
are written as
\begin{align}
  \label{submatrix_eq-Both_Dirichlet}
  \left(
   \begin{array}{cccc}
      -g_{\NodeAdomI}^{\NodeAdomI}
     &-g_{\NodeAdomI}^{\NodeBdomI}
     &0
     &h_{\NodeAdomI}^{\NodeBdomI}
     \\
      -g_{\NodeBdomI}^{\NodeAdomI}
     &-g_{\NodeBdomI}^{\NodeBdomI}
     &0
     &c_{\NodeBdomI}^{}+h_{\NodeBdomI}^{\NodeBdomI}
     \\
      0
     &\kappa_{21} g_{\NodeAdomII}^{\NodeBdomII}
     &-g_{\NodeAdomII}^{\NodeAdomII}
     &h_{\NodeAdomII}^{\NodeBdomI}
     \\
       0
     &\kappa_{21} g_{\NodeBdomII}^{\NodeBdomII}
     &-g_{\NodeBdomII}^{\NodeAdomII}
     &c_{\NodeBdomII}^{}+h_{\NodeBdomII}^{\NodeBdomI}
   \end{array}
  \right)
  \left(
   \begin{array}{c}
      q_{\NodeAdomI}^{}\\q_{\NodeBdomI}^{}\\q_{\NodeAdomII}^{}\\u_{\NodeBdomI}^{}
   \end{array}
  \right).
\end{align}
By applying \EQthree{c-ia-ib}{h-ia-ib-sing}{g-ia-ib} to the sub-matrix
to eliminate coefficients with the subscript $\NodeB^{(1)}$ and $\NodeB^{(2)}$,
the sub-matrix is rewritten as
\begin{align}
  \label{submat1-Dirichlet-case}
  \left(
   \begin{array}{cccc}
       -g_{\NodeAdomI}^{\NodeAdomI}
     &-g_{\NodeAdomI}^{\NodeBdomI}
     &0
     &0
     \\
      -g_{\NodeAdomI}^{\NodeAdomI}
     &-g_{\NodeAdomI}^{\NodeBdomI}
     &0
     &c_{\NodeAdomI}^{}
     \\
      0
     &\kappa_{21} g_{\NodeAdomII}^{\NodeBdomII}
     &-g_{\NodeAdomII}^{\NodeAdomII}
     &0
     \\
      0
     &\kappa_{21} g_{\NodeAdomII}^{\NodeBdomII}
     &-g_{\NodeAdomII}^{\NodeAdomII}
     &c_{\NodeAdomII}^{}
   \end{array}
  \right).
\end{align}
Furthermore, by applying elementary matrix operations,
the sub-matrix is transformed to
\begin{align}
  %\to
  \left(
   \begin{array}{cccc}
      -g_{\NodeAdomI}^{\NodeAdomI}
     &-g_{\NodeAdomI}^{\NodeBdomI}
     &0
     &0
     \\
      0
     &0
     &0
     &c_{\NodeAdomI}^{}
     \\
       0
     &\kappa_{21} g_{\NodeAdomII}^{\NodeBdomII}
     &-g_{\NodeAdomII}^{\NodeAdomII}
     &0
     \\
      0
     &0
     &0
     &c_{\NodeAdomII}^{}
   \end{array}
  \right),
\end{align}
where the second row is the result of replacement
with the difference between the first and the second row of the matrix in \EQ{submat1-Dirichlet-case},
and the fourth row is obtained by similar operations.
The non-zero components in both the second and the fourth row are only found in the fourth column;
therefore, the rank of this sub-matrix is reduced by one.

In other cases,
the product of the sub-matrices and the unknown vector
after applying \EQthree{c-ia-ib}{h-ia-ib-sing}{g-ia-ib} are as follows.
In the case of $\NodeA^{(1)}\in I_N$ and $\NodeA^{(2)}\in I_N$,
\begin{align}
  \label{submatrix_eq-Both_Neumann}
  \left(
   \begin{array}{cccc}
      c_{\NodeAdomI}^{}
     &0
     &-g_{\NodeAdomI}^{\NodeBdomI}
     &0
     \\
      0
     &c_{\NodeAdomI}^{}
     &-g_{\NodeAdomI}^{\NodeBdomI}
     &0
     \\
      0
     &0
     &\kappa_{21} g_{\NodeAdomII}^{\NodeBdomII}
     &c_{\NodeAdomII}^{}
     \\
      0
     &c_{\NodeAdomII}^{}
     &\kappa_{21} g_{\NodeAdomII}^{\NodeBdomII}
     &0
   \end{array}
  \right)
  \left(
   \begin{array}{c}
      u_{\NodeAdomI}^{}\\u_{\NodeBdomI}^{}\\q_{\NodeBdomI}^{}\\u_{\NodeAdomII}^{}
   \end{array}
  \right).
\end{align}
In the case of $\NodeA^{(1)}\in I_D$ and $\NodeA^{(2)}\in N_N$,
\begin{align}
  \label{submatrix_eq-coupled_Dirichlet_Neumann}
  \left(
   \begin{array}{cccc}
       -g_{\NodeAdomI}^{\NodeAdomI}
     &-g_{\NodeAdomI}^{\NodeBdomI}
     &0
     &0
     \\
      -g_{\NodeAdomI}^{\NodeAdomI}
     &-g_{\NodeAdomI}^{\NodeBdomI}
     &c_{\NodeAdomI}^{}
     &0
     \\
      0
     &\kappa_{21} g_{\NodeAdomII}^{\NodeBdomII}
     &0
     &c_{\NodeAdomII}^{}
     \\
      0
     &\kappa_{21} g_{\NodeAdomII}^{\NodeBdomII}
     &c_{\NodeAdomII}^{}
     &0
   \end{array}
  \right)
  \left(
   \begin{array}{c}
      q_{\NodeAdomI}^{}\\q_{\NodeBdomI}^{}\\u_{\NodeBdomI}^{}\\u_{\NodeAdomII}^{}
   \end{array}
  \right).
\end{align}
Since $\kappa_{21}>0$ generally,
\EQtwo{submatrix_eq-Both_Neumann}{submatrix_eq-coupled_Dirichlet_Neumann} do not become singular.

Furthermore, in the case of $\NodeA^{(1)}\in I_R$ or $\NodeB^{(2)}\in I_R$,
the coefficient matrix becomes
more complex than \EQ{submatrix_eq-Both_Neumann} or \EQ{submatrix_eq-coupled_Dirichlet_Neumann};
therefore, the rank is also not reduced.

\subsubsection{Case of internal interface condition}
\label{sec:internal_interface-case-CBIE}
In the case where the node is located
at a corner of two
interface boundary elements
as shown in \Fig{double_node-interface-domain}(a),
there are four sub-nodes at the corner node.
After applying \EQthree{c-ia-ib}{h-ia-ib-sing}{g-ia-ib},
and eliminating the quantities
related to $\NodeA^{(2)}$ and $\NodeB^{(2)}$,
the terms related to the sub-nodes in the four sub-node equations are expressed as
\begin{align}
  \label{submatrix_eq-Interface}
  \left(
  \begin{array}{cccc}
              -g_{\NodeAdomI}^{\NodeAdomI}
    &
              -g_{\NodeAdomI}^{\NodeBdomI}
    &
    c_{\NodeAdomI}
    &
    0
    \\
              -g_{\NodeAdomI}^{\NodeAdomI}
    &
              -g_{\NodeAdomI}^{\NodeBdomI}
    &
    0
    &
    c_{\NodeAdomI}
    \\
    \kappa_{21}g_{\NodeAdomII}^{\NodeAdomII}
    &
    \kappa_{21}g_{\NodeAdomII}^{\NodeBdomII}
    &
    c_{\NodeAdomII}
    &
    0
    \\
    \kappa_{21}g_{\NodeAdomII}^{\NodeAdomII}
    &
    \kappa_{21}g_{\NodeAdomII}^{\NodeBdomII}
    &
    0
    &
    c_{\NodeAdomII}
  \end{array}\right)
  \left(
  \begin{array}{cccc}
    q_{\NodeAdomI}\\q_{\NodeBdomI}\\u_{\NodeAdomI}\\u_{\NodeBdomI}
  \end{array}\right).
\end{align}
By applying the elementary matrix operations to the coefficient matrix,
the matrix is transformed to
\begin{align}
  \left(
  \begin{array}{cccc}
    -           g_{\NodeAdomI}^{\NodeAdomI}
    &
    -           g_{\NodeAdomI}^{\NodeBdomI}
    &
    c_{\NodeAdomI}
    &
    0
    \\
    0
    &
    0
    &
    c_{\NodeAdomI}
    &
    -c_{\NodeAdomI}
    \\
    \kappa_{21}g_{\NodeAdomII}^{\NodeAdomII}
    &
    \kappa_{21}g_{\NodeAdomII}^{\NodeBdomII}
    &
    c_{\NodeAdomII}
    &
    0
    \\
    0
    &
    0
    &
    c_{\NodeAdomII}
    &
    -c_{\NodeAdomII}
  \end{array}\right).
\end{align}
We can understand that the rank of this matrix is reduced by one
since the determinant of a 2x2 sub-matrix composed of non-zero elements
in the second and the fourth row in the above sub-matrix is zero.
When we increase the number of the domains as shown in \Fig{double_node-interface-domain}(b),
we can illustrate the singularity, although the details are not shown here.
Therefore, the sub-node equations become singular when the node
is located at an intersection of the interface boundaries
which partitions the original single domain into two or more domains.

%%%%%%

\subsection{Summary of rank deficiency conditions for CBIEs}
\label{sec:CBIE-summary}
In this section, we analyze the rank deficiency conditions for the sub-node equations of a double node for CBIEs.
The rank of the coefficient matrix is reduced in the following cases:
\begin{itemize}
  \item both sub-nodes of a double node are imposed by Dirichlet conditions,
  \item sub-nodes are located at an intersection of
      an interface boundary element
      and two boundary elements imposed by Dirichlet conditions,
  \item sub-nodes are located at intersections
        of interface boundary elements with different normal directions
        and they do not belong to any boundary with ordinary boundary conditions.
\end{itemize}
The original sub-node equations or the sub-matrix associated with the sub-nodes for these cases
include
a pair of \EQtwo{Both-Dirichlet-ia}{Both-Dirichlet-ib} for the first case,
\EQ{submatrix_eq-Both_Dirichlet} for the second case,
and \EQ{submatrix_eq-Interface} for the last case.
These equations have common characteristics;
all the equations include $q_\NodeA$ and $q_\NodeB$ (or $q_\NodeAdomI$ and $q_\NodeBdomI$ for the interface boundary) as the unknown variables.
The coefficients associated with these unknowns are $g_i^j$ ($i,j\in\{\NodeA,\NodeB\}$),
which possess the characteristics shown in \EQ{g-ia-ib}.
The terms with $q_\NodeA$ and $q_\NodeB$
are canceled because of this characteristic during the transformations of the equations,
and the ranks have been reduced.

If the coefficients $g_i^j$ did not satisfy \EQ{g-ia-ib};
in other words, if they
included information on $\Vn_i$, these singularities could be avoided.

%%%%%%%%%%%%%%%%%%%%%%%%%%%%%%%%%%%%%%%%%%%%%%%%%%%%%%%%%%%%%%%%

\section{Regularization of coefficient matrix using HBIEs}
\label{sec:HBIE}
As described in Sec.~\ref{sec:CBIE-summary},
the reason for a rank deficient matrix in CBIE is
the coefficients $g_i^j$ do not carry information concerning the normal direction of the node $i$.
To provide this information,
we employ a hypersingular boundary integral equation (HBIE).
The HBIE can be derived by taking a gradient of the CBIE shown in \EQ{bie1}
with respect to $\Vr_i$:
\begin{align}
  \label{grad-bie1}
    \Gradxx{\left[c(\xx)u(\xx)\right]}
    =\GOint{}{&\left[(\Gradxx{ u^*(\x,\xx)})(\Grad{u(\x)})\IP\Vn
     \right.\nonumber\\&\left.
              -u(\x) (\Gradxx{\Grad{ u^*(\x,\xx)}})\IP\Vn \right]},
    %\nonumber\\
\end{align}
where $\nablaxx$ denotes the gradient with respect to $\xx$.
The gradient, $\Gradxx\Grad u^*$, has a stronger singularity
than the CBIE,
and this is known as hypersingularity.
However, the integral of the hypersingular function can be regularized.
In this study,
we employ the regularization
for the Helmholtz equation which uses the fundamental solution of Laplace's equation \cite{Tomioka:2010}.

\subsection{Regularized gradient field}
In this section,
after presenting a regularization of the gradient at $\xx$
according to the method shown in Ref.~\cite{Tomioka:2010},
a discretized form of the normal derivative, $q$, is derived.

The gradient is expressed by a 2x2 dyadic tensor $\DCi$ and a vector $\Vector{J}_i(u,q)$ as
\begin{align}
  \label{Reg-grad-bie}
   &
  \DCi\,\IP\Grad{u(\xx)}+\Vector{J}_i(u,q)=\Vector{0}.
\end{align}
The second term can be evaluated by two types of boundary integrals as
\begin{align}
  &
  \label{J_i-def}
  \Vector{J}_i(u,q)=
      \sum_{\gamma=\EleA,\EleB}\Vector{J}_i^{\gamma,\rm sing}(u,q)
     +\sum_{k\ne \EleA,\EleB}\Vector{J}_{i}^{k,\rm reg}(u,q),
\end{align}
where $\Vector{J}_i^{\gamma, \rm sing}(u,q)$ is associated to the boundary integral
along the boundary elements that include the node $i\in\{\NodeA,\NodeB\}$
and $\gamma\in\{\EleA,\EleB\}$ (\Fig{double_node});
whereas $\Vector{J}_i^{k,\rm reg}$ corresponds to integrals
of the boundary elements not including the node $i$.
Similar to the factor $c(\xx)$ in CBIEs,
the singular integral included in \EQ{grad-bie1}
contributes to $\DCi$ on the left-hand side of \EQ{Reg-grad-bie}
and $\Vector{J}_i^{\gamma,\rm sing}(u,q)$.
The components of $\Vector{J}_i^{\gamma,\rm sing}(u,q)$, $\Vector{J}_{i}^{k,\rm reg}(u,q)$ and $\DCi$
are written as
\begin{align}
  \label{J-AB-def}
  &\Vector{J}_i^{\gamma,\rm sing}(u,q)
   =
    -\Rint{0}{L_\gamma}{\frac1r\left(\RoundR{ u^*}-\RoundR{u_L^*}\right)}\cdot u(\xx)\Vnal 
  \nonumber\\&\hspace*{3em}
    +\left.\Rint{0}{L_\gamma}{r\RoundR{ u^*}}\right.
     \left(
      -\frac12\Vnal\FuncAt{\ROUND{^2u}{\sal^2}}{\xx}
      +\Vsal\FuncAt{\frac{\partial^2u}{\partial\sal\partial\nal}}{\xx}
     \right),
  %%%%%%
  \\
  \label{Jn-def}
  &\Vector{J}_i^{k,\rm reg}(u,q)
   =
       \Gint{k}{\left\{(\Gradxx q^*)u-(\Gradxx q_L^*) u(\xx) -(\Gradxx u^*)q \right\}},
  \\
  \label{DC-def}
  &\DCi=
    %\nonumber\\
    %&
    \left(\begin{array}{cc}
       \displaystyle
       \Diff{\textstyle\frac{2\theta_\gamma+\sin{2\theta_\gamma}}{4\pi}}
        &
       \hspace*{-1.0em}%
       \displaystyle
       \Diff{\textstyle\frac{-\cos{2\theta_\gamma}}{4\pi}+ u^*(L_\gamma)}
       \\[1em]
       \displaystyle
       \Diff{\textstyle\frac{-\cos{2\theta_\gamma}}{4\pi}- u^*(L_\gamma)}
        &
       \hspace*{-1.0em}%
       \displaystyle
       \Diff{\textstyle\frac{2\theta_\gamma-\sin{2\theta_\gamma}}{4\pi}}
    \end{array}\right)
    ,
   \\
  &\Diff{f_\gamma}\Define f_{\EleA}-f_{\EleB},
\end{align}
where $u_L^*$ denotes the fundamental solution to the Laplace equation that satisfies \EQ{funda-laplace},
$q_L^*$ is its normal derivative with respect to $\Vn(\x)$ at $\x$;
and $\Vector{\tau}_\gamma$, $\theta_\gamma$,
and $L_\gamma$ are the tangential unit vectors along $\Gamma_\gamma$,
the azimuth angle of $\Gamma_\gamma$ from $\xx$,
and the length of $\Gamma_\gamma$, respectively.

Similar to CBIEs,
the non-singular integral in \EQ{Jn-def} is expressed by the shape function as
\begin{align}
   \label{s_int-reg}
     &
     \Vector{s}^{(k,l),\rm reg}_{i}=\Gint{k}{\nabla_i\nabla u^*(\x,\xx)\phi^{(k,l)}(\x)\IP\Vn},
   \\
   \label{sL_int-reg}
     &
     \Vector{s}^{(k,l),\rm reg}_{L,i}=\Gint{k}{\nabla_i\nabla u_L^*(\x,\xx)\phi^{(k,l)}(\x)\IP\Vn},
   \\
   \label{t_int-reg}
     &
     \Vector{t}^{(k,l),\rm reg}_{i}=\Gint{k}{\nabla_i u^*(\x,\xx)\phi^{(k,l)}(\x)},
\end{align}
which are evaluated using numerical integrals.
By using the permutation matrix,
the second term on the left-hand side in \EQ{J_i-def}
is denoted by
a form of linear combinations as
\begin{align}
   \label{J-reg-lin}
   \sum_{k\ne \EleA,\EleB}\Vector{J}_i^{k,\rm reg}(u,q)
     &
     =\sum_{j\ne \NodeA,\NodeB}\Vector{s}_i^{j,\rm reg}u_j^{}
     -\sum_{j\ne \NodeA,\NodeB}\Vector{s}_{L,i}^{j,\rm reg}u_i^{}
     -\sum_{j\ne \NodeA,\NodeB}\Vector{t}_i^{j,\rm reg}q_j^{}.
\end{align}

The two singular integrals on the right-hand side of \EQ{J-AB-def} can be analytically obtained
by using Taylor series expansions around $\xx$ for the two second-order derivatives.
The first term on the right-hand side of \EQ{J_i-def}
is also expressed as linear combinations with respect to the nodes neighboring the sub-nodes $\NodeA$ and $\NodeB$
as
\begin{align}
   \label{J-sing-lin}
     &
   \sum_{\gamma=\EleA,\EleB}\Vector{J}_i^{\gamma,\rm sing}(u,q)
     =\sum_{j\in \tilde{I}^{\rm sing}}\Vector{s}_i^{j,\rm sing}u_j^{}
     -\sum_{j\in \tilde{I}^{\rm sing}}\Vector{t}_i^{j,\rm sing}q_j^{},
     \nonumber\\
    &\qquad
     \tilde{I}^{\rm sing}=\{i|\xx\in\Gamma_\EleA\cup\Gamma_\EleB\}=\{(\EleA,1),(\EleA,2),(\EleB,1),(\EleB,2)\}.
\end{align}
By defining
\begin{align}
   \Vector{s}_i^j=\Vector{s}_i^{j,\rm sing}+\Vector{s}_i^{j,\rm reg}-\Vector{s}_{L,i}^{j,\rm reg},
   \quad
   \Vector{t}_i^j=\Vector{t}_i^{j,\rm sing}+\Vector{t}_i^{j,\rm reg},
\end{align}
the vector $\Vector{J}_i(u,q)$ is expressed as the following form with a linear combination:
\begin{align}
  \label{J-sum_form}
  \Vector{J}_i(u,q)=\sum_j\Vector{s}_i^{j}u_j-\sum_j\Vector{t}_i^{j}q_j^{}.
\end{align}
By substituting this equation into \EQ{Reg-grad-bie} and by multiplying both sides by $\DCiinv$,
the gradient field at $\xx$ is expressed as
\begin{align}
  \label{grad-by-C}
  \Grad{u(\xx)}=-\DCiinv\,\IP\left(\sum_j\Vector{s}_i^{j}u_j-\sum_j\Vector{t}_i^{j}q_j^{}\right).
\end{align}

\subsection{Rank deficiency conditions in HBIEs}
\label{sec:Rank_deficient_condition-HBIE}
Let us consider the relations of $\Vector{s}_i^j$ and $\Vector{t}_i^j$ when $\xx$ is a double node;
i.e., $\x_i=\x_\NodeA=\x_\NodeB$.
If both $u_j$ and $q_j$ at all nodes are determined as known values,
the gradient in \EQ{grad-by-C} must be expressed uniquely.
Since the dyadic tensor $\DCiinv$ on the right-hand side of \EQ{grad-by-C}
is determined by a geometric configuration of the boundary elements connected
to $\xx$ as shown in \EQ{DC-def},
$\Dyad{C_{\mathit\NodeA}}=\Dyad{C_{\mathit\NodeB}}$ is satisfied.
Therefore, to uniquely determine the gradient, i.e., $\Grad{u(\x_\NodeA)}=\Grad{u(\x_\NodeB)}$,
the coefficient vectors for $\NodeA$ and $\NodeB$ must be the same.
Strictly speaking, there is an exception in $\Vector{s}_i^j$ for $j=\NodeA,\NodeB$.
In this case, the condition $\Vector{s}_{\NodeA}^j\ne\Vector{s}_{\NodeB}^j$
is permitted
when $\Vector{s}_{\NodeA}^{\NodeA}+\Vector{s}_{\NodeA}^{\NodeB}=\Vector{s}_{\NodeB}^{\NodeA}+\Vector{s}_{\NodeB}^{\NodeB}$ is satisfied
since $u_{\NodeA}=u_{\NodeB}$.
Therefore, 
\begin{align}
  \label{s-cond-sing}
  &\Vector{s}_{\NodeA}^j\ne\Vector{s}_{\NodeB}^j && \mbox{ for $j\in\{\NodeA,\NodeB\}$},\\
  \label{s-cond}
  &\Vector{s}_{\NodeA}^j=\Vector{s}_{\NodeB}^j && \mbox{ for $j\not\in\{\NodeA,\NodeB\}$},\\
  \label{t-cond}
  &\Vector{t}_{\NodeA}^j=\Vector{t}_{\NodeB}^j && \mbox{ for any $j$}.
\end{align}

Since $q_i=\Vn(\xx)\cdot\Grad{u(\xx)}$,
by taking an inner product of $\Vn(\xx)$ to \EQ{grad-by-C}, the discretized equation for the HBIE is obtained as
\begin{align}
 \label{HBIE-point_i}
  &
  \sum_j v_i^j u_j^{} -\sum_j(-\delta_i^j+w_i^j)q_j^{}=0,
\end{align}
where
\begin{align}
 \label{def-v_i^j}
  &
  v_i^j=\Vn_i\,\IP\,\DCiinv\,\IP\,\Vector{s}_i^j,
  \\
 \label{def-w_i^j}
  &
  w_i^j=\Vn_i\,\IP\,\DCiinv\,\IP\,\Vector{t}_i^j.
\end{align}

When $\Vn_\NodeA\ne\Vn_\NodeB$,
we can derive the following conditions using \EQtwo{s-cond}{t-cond}:
\begin{align}
  \label{v-cond}
  & v_{\NodeA}^j\ne v_{\NodeB}^j && \mbox{ when $j\ne\NodeA,\NodeB$},\\
  \label{w-cond}
  & w_{\NodeA}^j\ne w_{\NodeB}^j && \mbox{ for any $j$}.
\end{align}
These relations mean that
the two sub-node equations of \EQ{HBIE-point_i} for $i=\NodeA$ and $i=\NodeB$
are different.
Even in the case of $\Vn_\NodeA=\Vn_\NodeB$, the two sub-node equations are also different
since
the coefficients for $u_\NodeA$ of the two equations are always different because of \EQ{s-cond-sing} and the coefficients for $q_\NodeA$
are always different because of the nature of $\delta_i^j$.

The aforementioned difference in the two equations is not a sufficient condition to not reduce the rank
since we have not proved that
the determinant of the sub-matrix associated with the sub-nodes is not equal to zero for any problem.
However, the HBIE is highly applicable to the problem
in which the node equations constructed by the CBIE become singular.

\begin{figure}[tb]
  \centering
  \includegraphics[scale=0.95]{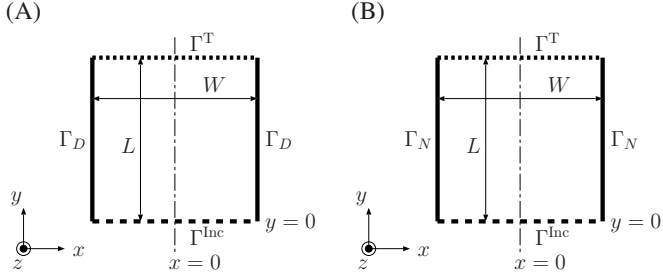}
  \caption{Analysis models.
    (A) Dirichlet conditions of $\overline{u}=0$ on $\Gamma_D$
        and an incident wave condition with a mode proportional to $\cos(k_x x)$ on $\Gamma^{\rm Inc}$ are given.
    (B) Neumann conditions of $\overline{q}=0$ on $\Gamma_N$
        and an incident condition with $\sin(k_x x)$ on $\Gamma^{\rm Inc}$.
    In both cases, either of the three termination conditions on $\Gamma^{\rm T}$ is given;
    (a) a short condition, (b) an open condition, or (c) a non-reflective condition.
  }
  \label{fig:model}
\end{figure}

\begin{figure}[tb]
  \centering
  \includegraphics[width=0.97\hsize]{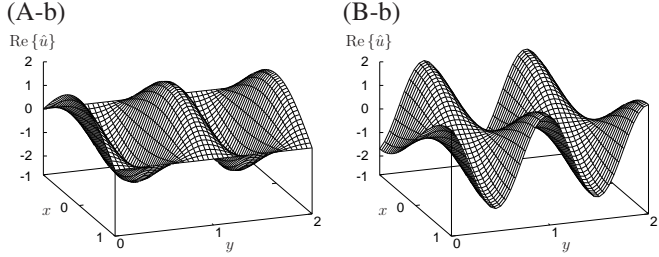}
  \caption{Exact solutions for model (A-b) and model (B-b).
    In both models, the termination condition at $y=2$ is $R=+1$ (conditions for the open termination:(b)).
    The conditions $x=\pm1$ are Dirichlet conditions with $\overline{u}=0$ for (A)
    and Neumann conditions with $\overline{q}=0$ for (B).
    Although the solution is a complex valued function, this figure only illustrates the real part.
    The imaginary part obeys the same eigenfunction, but the amplitude is different from the real part.
  }
  \label{fig:exact_dist-open_end-1reg}
\end{figure}

\begin{figure}[tb]
  \centering
  \includegraphics[width=0.97\hsize]{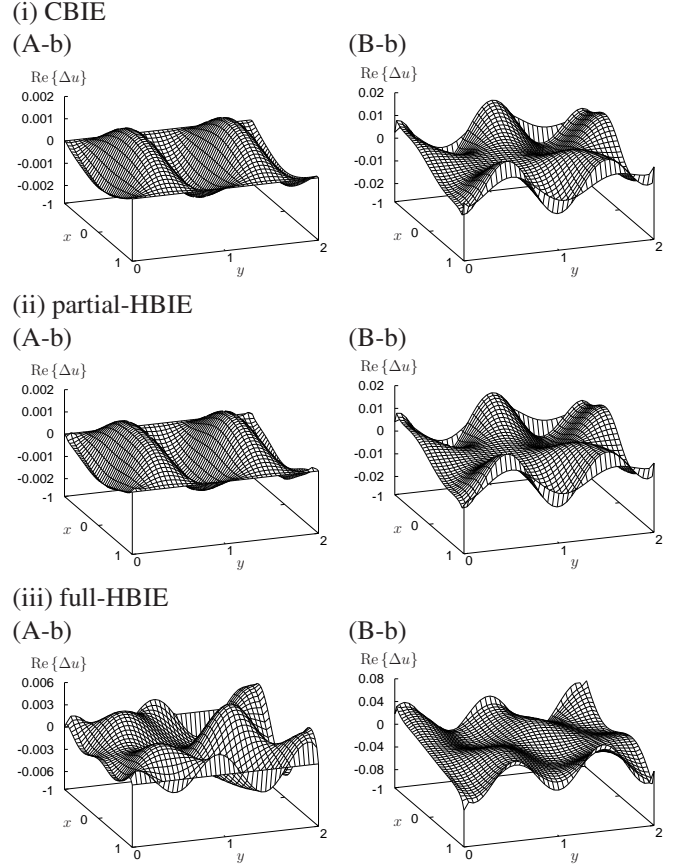}
  \caption{Error distributions of numerical solutions:
   (i) CBIE, (ii) partial-HBIE, and (iii) full-HBIE.
   The boundary condition on $\Gamma^{\rm T}$ is an open termination condition (model (b));
   and on $x=\pm1$, Dirichlet conditions (A), and Neumann conditions (B).
   The size of the boundary elements is 0.05.
   Each sub-figure illustrates the real part of the error, ${\rm Re}\left\{\Delta u\right\}$.
   The scales of the vertical axes are different between sub-figures.
  }
  \label{fig:error_dist-open_end-1reg}
\end{figure}

\begin{figure}[tbh]
  \centering
  \includegraphics[width=0.7\hsize]{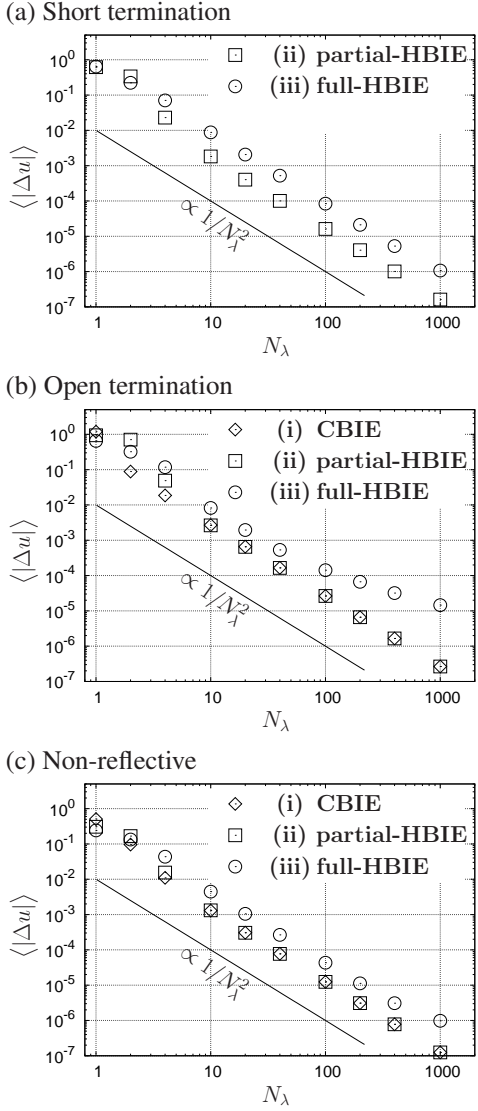}
  \caption{%
    Element size dependence of error.
    (a) Short termination, (b) Open termination, and (c) Non-reflective condition.
    The horizontal axis, $N_\lambda$, denotes the number of elements in a wavelength;
    i.e., $N_\lambda=\lambda/\Delta_{\rm Ele}$
    where $\Delta_{\rm Ele}$ is the boundary element size.
    In sub-figure (a), the results by CBIE (i) are not shown
    as the analyses failed because of the singularity of the coefficient matrix.
  }
  \label{fig:size-vs-error}
\end{figure}

\begin{figure}[tb]
  \centering
  \includegraphics[width=0.7\hsize]{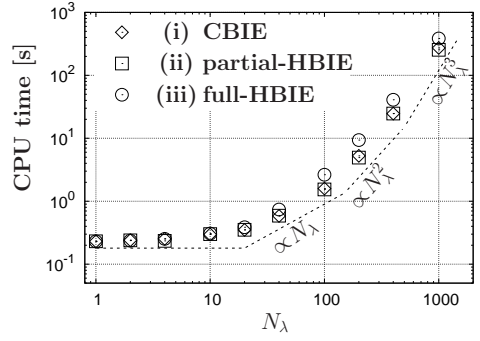}%
  \caption{%
    Element size dependence of CPU time.
    The horizontal axis, $N_\lambda$, denotes the number of elements in a wavelength. 
    The model is (A-b); i.e., Dirichlet conditions for the side boundaries,
    and open termination conditions for the top boundary.
    The CPU time for the other termination condition (a) and (c) are almost the same as this result;
    however, the result of the CBIE (i) is not obtained
    for the condition of the short termination (a).
  }
  \label{fig:size-vs-cpu}
\end{figure}

\begin{figure}[tb]
  \centering
  \includegraphics[width=0.99\hsize]{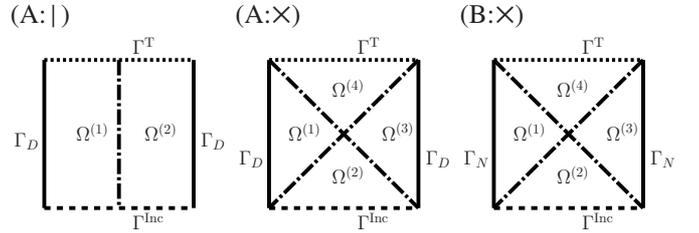}%
  \caption{%
    Analysis models including interfaces between media.
    The first letter in the label of each sub-figure presents the base model shown in \Fig{model},
    and the last symbol illustrates a shape of the interface boundaries drawn by the dot-dash lines.
    The media constants of each domain are the same.
  }
  \label{fig:model_with_4reg}
\end{figure}

\begin{figure}[tb]
  \centering
  \includegraphics[width=0.9\hsize]{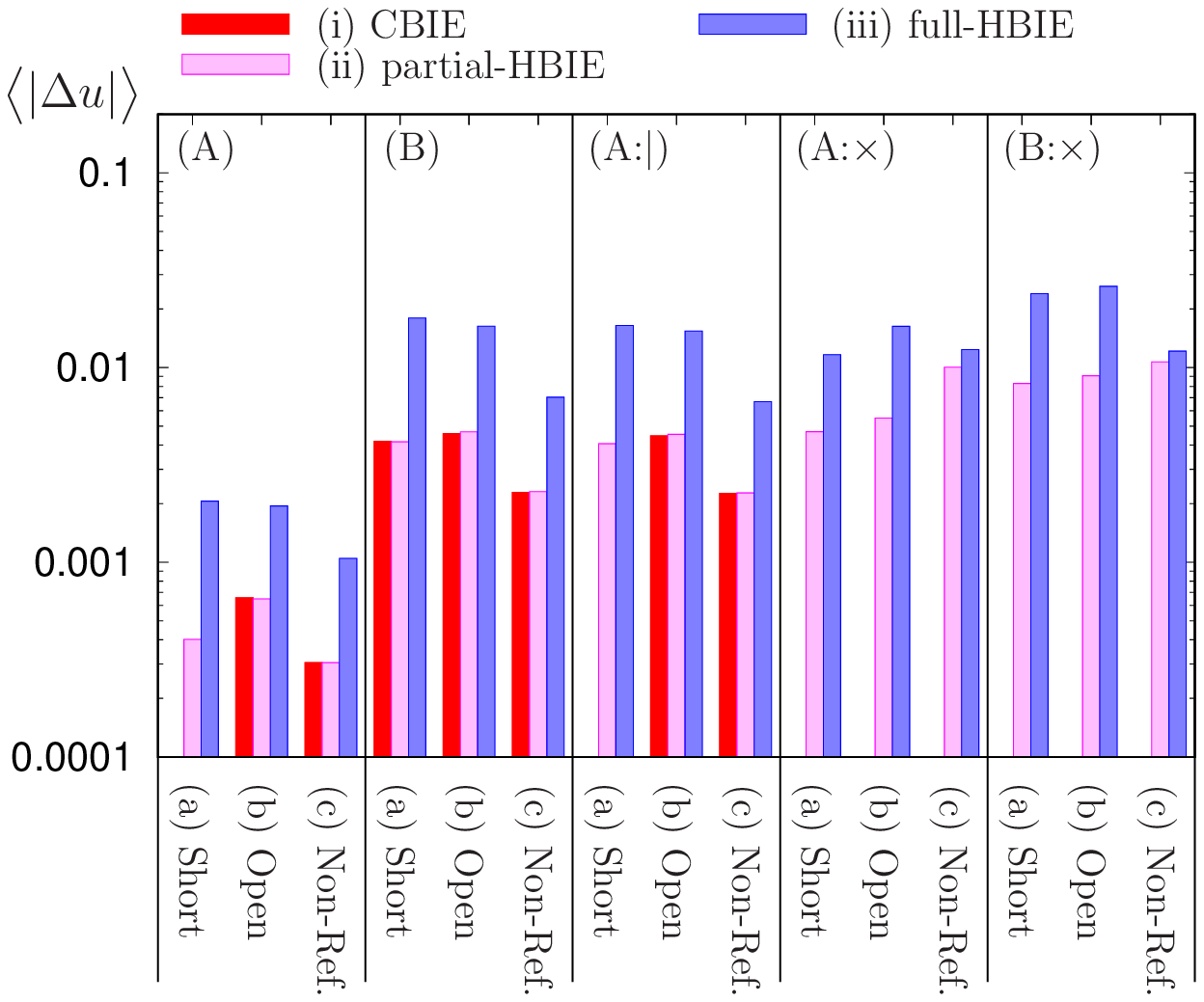}
  \caption{%
    Comparisons of errors in the problems including the interface boundaries.
    The label at the top of each block presents the model shown in \Fig{model} and \Fig{model_with_4reg}.
    In each model,
    the error obtained by the three methods, (i) CBIE, (ii) partial-HBIE, and (iii) full-HBIE,
    are expressed as a triplet of band graphs.
    The case in which the band graph is not shown means that the set of simultaneous equations is singular.
  }
  \label{fig:error-interface_conditions}
\end{figure}

%%%%%%%%%%%%%%%%%%%%%%%%%%%%%%%%%%%%%%%%%%%%%%%%%%%%%

\section{Numerical results and discussions}
\label{sec:results}
To demonstrate the validity of the method in regularizing the coefficient matrix using HBIEs,
we analyzed simple models as shown in \Fig{model}.
The actual models correspond to electromagnetic wave propagation problems in a waveguide
where the wall parallel to the propagation direction is made of metal.
When we consider $u$ as the $z$-component of the electric field,
$u=0$ at the side-walls of the waveguide
since the component of the electric field parallel to the metal is zero.
In the case of model-(A), the metal walls are located at $\Gamma_D$,
which is called the TE${}_{10}$ mode for a waveguide with rectangular cross-section and width $W$.
In the case of model-(B),
virtual walls $\Gamma_N$ with $q=0$ are located at $x=\pm W/2$,
which are equivalent to the placement of physical metal walls with $u=0$ at $x=\pm W$,
which is called the TE${}_{20}$ mode.
These incident conditions are given as functions of $x$ as 
\begin{align}
  & \overline{u^{\rm inc}}(x)=\overline{u^{\rm inc}_0}\cos(k_x x) &&\mbox{for model-(A)},
  \\
  & \overline{u^{\rm inc}}(x)=\overline{u^{\rm inc}_0}\sin(k_x x) &&\mbox{for model-(B)},
  \\
  & k_x=\frac{\pi}{W}.
\end{align}
In these modes, a propagation field to $\pm y$-directions are
proportional to $\exp({\mp{\rm j}k_y})$, respectively,
where 
\begin{align}
  & k_y^2=k^2-k_x^2.
\end{align}
By using this characteristic,
the incident boundary condition on $\Gamma^{\rm inc}$ can be rewritten as a Robin condition:
\begin{align}
  \label{incident-condition}
  {\rm j}k_y u+q=2{\rm j}k_y\overline{u^{\rm inc}}.
\end{align}
On the boundary $\Gamma^{\rm T}$ at $y=L$, we analyzed the following three conditions
for termination:
\begin{center}
\begin{tabular}{cllll}
  (a) &Short:& $u=0$ & ($R=-1$),\\
  (b) &Open:& $q=0$ & ($R=+1$),\\
  (c) &\multicolumn{1}{l}{Non-reflective:}
    &${\rm j}k_y u+q=0$ & ($R=0$),
\end{tabular}
\end{center}
where $R$ denotes the reflection coefficient of the electric field on $\Gamma_T$; and
the non-reflective condition in (c) is equivalent to the incident condition
with $\overline{u^{\rm inc}}=0$ in \EQ{incident-condition}.
Among the above three, termination-(a) (short type)
has two double nodes at $(\pm W/2,L)$,
where both sub-nodes have two Dirichlet conditions;
therefore, the rank of the coefficient matrix of CBIEs will be reduced
as mentioned in Sec.~\ref{sec:Dirichlet-case-CBIE}.
The exact solutions to these models are given by
\begin{align}
  \label{exact-TE10}
  &
  \hat{u}(x,y)=\overline{u^{\rm inc}_0}\cos(k_x x)
    \left(e^{-{\rm j}k_y y}+R e^{{\rm j}k (y-2L)}\right)
   && \mbox{for model-(A)},
  \\
  \label{exact-TE20-central_half}
  &
  \hat{u}(x,y)=\overline{u^{\rm inc}_0}\sin(k_x x)
    \left(e^{-{\rm j}k_y y}+R e^{{\rm j}k (y-2L)}\right)
   && \mbox{for model-(B)}.
\end{align}
To simplify, we used $\overline{u^{\rm inc}_0}=1$, $W=2$, $L=2$, and $\lambda=2\pi/k=1$;
under this condition,
$\max|u|=2$ for the cases (a) and (b), and $\max|u|=1$ for the case (c),
$\lambda_x=2\pi/k_x=4$
and $\lambda_y=2\pi/k_y\simeq 1.033$.
As examples of $\hat{u}(x,y)$, the real part of $\hat{u}(x,y)$ for
the models (A-b) and (B-b) are shown in \Fig{exact_dist-open_end-1reg}.

We analyzed three types of the following simultaneous equations:
\begin{enumerate}
  \def\labelenumi{(\roman{enumi})}
  \item CBIE:
     BIEs for all nodes including sub-nodes are obtained from the CBIE,
  \item partial-HBIE:
     BIEs for only sub-nodes related to double nodes are obtained from the HBIE,
     and BIEs for the other nodes are obtained from the CBIE.
  \item full-HBIE:
     BIEs for all nodes including sub-nodes are obtained from the HBIE,
\end{enumerate}
To solve the complex-valued simultaneous equations,
we employed two subroutines based on an LU decomposition provided in Lapack \cite{Lapack},
in which the subroutine names are `zgetrf' and `zgetrs.'
After determining both $u$ and $q$ at all boundary nodes
by solving the simultaneous equations,
the internal field $u(\x_{i\,'})$ can be evaluated by 
\begin{align}
  \label{internal-field}
  u(\x_{i\,'})=\sum_{j\in I} g_{i\,'\,}^j q_j -\sum_{j\in I}h_{i\,'\,}^j u_j,
\end{align}
which is derived from the discretized CBIE shown in \EQ{bie1}.
This process is the same for all types of simultaneous equations. 

Figure \ref{fig:error_dist-open_end-1reg} presents the
distributions of errors
based on the three types of simultaneous equations (i), (ii), and (iii)
for the model-(A-b) and model-(B-b).
In these models with termination-(b) (open type),
the coefficient matrices, even in the case of the CBIE shown in \Fig{error_dist-open_end-1reg}(i),
are not singular as mentioned in Sec.~\ref{sec:CBIE-summary}.

Based on the comparison between the sub-figures (A-b-i) and (A-b-iii)
or between (B-b-i) and (B-b-iii) in \Fig{error_dist-open_end-1reg},
we can observe that the error of the full-HBIE is several times larger than that the CBIE.
This difference can be explained by two reasons.
The first reason is the difference in singularity to evaluate the coefficients.
The coefficients $w_i^j$ and $v_i^j$ in the HBIE
are evaluated from $\Vector{s}_i^j$ and $\Vector{t}_i^j$ with a multiplication of dyadic $\DCi$
as shown in \EQtwo{def-v_i^j}{def-w_i^j}, respectively.
The vectors $\Vector{s}_i^j$ and $\Vector{t}_i^j$ are evaluated
by the boundary integral shown in \EQthree{s_int-reg}{sL_int-reg}{t_int-reg}.
Their integrands include the first or second order derivatives of the fundamental solution.
The strongest singularity in the HBIE is $O(r^{-2})$,
while the strongest singularity in the CBIE is $O(r^{-1})$.
This error emerges significantly in the contributions between two nodes with short distances.
The second reason is the multiplication by $\DCi$.
In the worst case, the errors of $v_i^j$ and $w_i^j$ are multiplied by the maximum norm, $\Norm{\DCiinv}$, and the errors of $\Vector{s}_i^j$ and $\Vector{t}_i^j$, respectively.
This amplification affects all coefficients regardless of the distances between nodes.
According to Ref.~\cite{Tomioka:2010},
$\Norm{\DCiinv}\le 4\pi/(\pi-2)\simeq 11$ in the case of $\Delta\theta=\pi/2$ and $L_{\EleA}=L_{\EleB}$.

In contrast,
the error in the case of the partial-HBIE shown in \Fig{error_dist-open_end-1reg}(A-b-ii)
is similar to the case of the CBIE in (A-b-i);
and the relation between (B-b-ii) and (B-b-i) is also similar.
The number of unknowns 
in the analyses presented in \Fig{error_dist-open_end-1reg} was 164 including four double nodes;
i.e., the number of HBIEs was only 8 and the number of CBIEs was 156 in the partial-HBIE.
Since the number of HBIEs with a larger error is sufficiently smaller than that of CBIEs,
the total error of partial-HBIEs does not increase so much as that of CBIEs.

Based on the comparisons between model-(A) and model-(B) in \Fig{error_dist-open_end-1reg},
we can observe that the error of model-(B) is larger than that of model-(A) for each equation type.
This reason is the same as the first reason for the difference between CBIEs and full-HBIEs.
The coefficients $g_{i\,'}^j$ and $h_{i\,'}^j$ in \EQ{internal-field}
are results of the boundary integrals in which integrands include $u^*$ and $q^*$, respectively.
Since the singularity of $q^*$ is stronger than that of $u^*$,
the contribution of $h_{i\,'}^j$ to $u(\x_{i\,'})$ is larger than that of $g_{i\,'}^j$,
especially in the case
where the distance between the field point and boundary elements is smaller than a wavelength.
The error caused by $h_{i\,'}^j$ is also larger than $g_{i\,'}^j$.
In model-(A), there are many nodes with $\overline{u}=0$.
Therefore, the larger error caused by $h_{i\,'}^j$ does not appear,
and the smaller error caused by $g_{i\,'}^j$ becomes dominant.
In contrast, in model-(B),
since there are no nodes with $\overline{u}=0$,
the larger error caused by $h_{i\,'}^j$ remains.

We analyzed the errors of different boundary element sizes.
The error is evaluated based on average sampling points as follows:
\begin{align}
  &
  \langle|\Delta u|\rangle=\frac{1}{N_s}\sum_{i'=1}^{N_s}{|\Delta u_{i'}|},
  \qquad
  \Delta u_{i'}=u_{i'}-\hat{u}_{i'},
\end{align}
where $N_s$ is the number of sampling points
that are intersections of the grid in \Fig{exact_dist-open_end-1reg} ($N_s=41^2$),
and $N_s$ is unchanged for all results regardless of the element size.
The errors are not normalized
since the averaged intensities of the exact solutions, $\langle|\hat{u}|\rangle$,
have almost the same order of magnitude;
$\simeq$0.8 for the termination types (a) and (b),
and $\simeq$0.6 for (c).
Figure~\ref{fig:size-vs-error} illustrates
the dependence between the error and the number of elements in a wavelength,
$N_\lambda=\lambda/\Delta_{\rm Ele}$ where $\Delta_{\rm Ele}$ denotes the boundary element size.
In the case of (a), there are no plots for the short termination type (i)
since the simultaneous equations become singular.
In all cases of (a), (b) and (c),
the errors decrease as $N_\lambda$ increases with a decay proportional to $1/N_\lambda^2$,
except for the points of $N_\lambda\gtrsim 200$ in (b).
This property is reasonable 
if $h_i^j$, $g_i^j$, $v_i^j$, and $w_i^j$ have accuracies of $O(\Delta_{\rm Ele})$
in the case of the linear element,
and the truncated error is proportional to $O(\Delta_{\rm Ele}^2)$.
However, there is an error which does not show this characteristic \cite{Tomioka:2010};
the order of the error of $\Vector{s}_i^{j,{\rm reg}}$ for a short distance
between the nodes $i$ and $j$ obeys $O(\Delta_{\rm Ele})$.
This error arises when $N_\lambda\gtrsim 200$.
In the cases (b) and (c),
the error of the full-HBIE (iii) is several times larger than the CBIE (i),
and that of the partial-HBIE (ii) is almost the same as (i),
which is similar to the result previously shown in \Fig{error_dist-open_end-1reg}.

Figure \ref{fig:size-vs-cpu} presents the computational time, which does not include the CPU time of the file input and output processes.
The computation consists of several major steps; 
the computation of the components of the coefficient matrix $a_i^j$
(through $h_i^j$ and $g_i^j$ in \EQtwo{h_int2}{g_int2}, respectively;
 or $v_i^j$ and $w_i^j$ in \EQtwo{def-v_i^j}{def-w_i^j}, respectively);
solving the matrix equation;
and the evaluation of the internal field using \EQ{internal-field};
for which individual computational costs are proportional
to $N_\lambda^2$, $N_\lambda^3$, and $N_\lambda$, respectively.
Higher-order terms appear with increasing $N_\lambda$.
In the case where $N_\lambda\lesssim5$,
most of the computation time is exhausted in minor common steps
such as initializing tables for the Hankel functions.
The computational cost for the full-HBIE (iii) is larger than that for the others.
This is because of the difference
between the evaluation time of $v_i^j$ and $w_i^j$ for the HBIE and that of $h_i^j$ and $g_i^j$ for the CBIE.
In the HBIE, the cost of evaluating $v_i^j$ and $w_i^j$ 
is mainly exhausted in the numerical integrals of non-singular elements for the three vectors
in \EQthree{s_int-reg}{sL_int-reg}{t_int-reg}.
In contrast, two scalar integrals in \EQtwo{h_int}{g_int} are dominant in the CBIE.
The cost of evaluating a coefficient with a vector is twice 
that of a scalar in two-dimensional problems,
and the number of components in the HBIE is 3/2 times greater than the CBIE.
Moreover,
the operator $\Vn\IP\Grad_i\Grad$ in \EQtwo{s_int-reg}{sL_int-reg} has
two vector components $\Ver\Ver\IP\Vn$ and $\Vn$.
Because some of the terms have common factors,
the sum of costs was reduced from these estimations;
however, the cost of evaluating the coefficient matrix component
in the HBIE is almost four times larger
than that in the CBIE.
Even when the simultaneous equations are singular,
it can be solved as a minimal-norm solution of underdetermined equations
by using a solver based on a singular value decomposition (SVD).
The details are not included in this paper because the authors do not understand
whether the minimal-norm solution is always correct or not.
By limiting the examples shown here, 
the accuracy of the CBIE by using a solver based on SVD
was almost the same as in the case of the partial-HBIE.
However, the computational time of SVD was much larger than in the case of LU decomposition;
e.g., the time to solve 6,900 s for the CBIE using SVD called `zgelss' in Lapack,
220 s for the partial-HBIE using the LU decomposition in the case of $N_\lambda=1,000$.

The above results can be summarized as follows.
First, as mentioned in Sec.~\ref{sec:Rank_deficient_condition-CBIE},
the set of simultaneous equations constructed by the CBIE for all nodes becomes singular
when both boundary conditions of the double node are imposed by Dirichlet conditions.
In contrast, when the node equations are constructed by the HBIE for all or a part of the nodes,
the set of simultaneous equations does not become singular.
Next, the accuracy of the HBIE is unfortunately worse than the CBIE when the set of equations is regular;
however, in the case of the partial-HBIE where
only the equations of the sub-nodes belonging to the double nodes
are given by the HBIE and the others are given by the CBIE,
the reduction in accuracy is negligibly small.
Finally,
more computational cost is required
to compute the component of the coefficient matrix by the HBIE
than the CBIE.
The rise in computational cost can be suppressed by applying the HBIE to sub-nodes only. 
Therefore, we can conclude that
the replacement of the CBIE by the HBIE only for the sub-nodes (partial-HBIE)
is the best solution from the viewpoint
of singularity, accuracy, and computational cost.

To demonstrate the applicability of the HBIE in the cases with interface boundaries,
we evaluated the models shown in \Fig{model_with_4reg}.
In these models,
the original model shown in \Fig{model}(A) or (B) is partitioned into
two or four sub-domains by interface boundaries,
and the exact solutions are the same as \EQtwo{exact-TE10}{exact-TE20-central_half}.
The number of sub-nodes for each double node is four at two intersections in the model-({A:\scalebox{1.0}{$\,|\,$}})
and at four corners in ({A:\scalebox{1.2}{$\times$}}) and ({B:\scalebox{1.2}{$\times$}}).
In addition, at the intersections of the four interface boundaries in ({A:\scalebox{1.2}{$\times$}}) and ({B:\scalebox{1.2}{$\times$}}), the number of sub-nodes for each double node is eight.
The combinations of four types of
boundary conditions (Dirichlet, Neumann, Robin, and interface conditions),
can be examined by these models.

Figure~\ref{fig:error-interface_conditions} presents comparisons of the errors.
As predicted in Sec.~\ref{sec:interface-case-CBIE},
when the simultaneous equations are constructed by the CBIE,
the set of equations for each model including the interface boundary conditions
is always singular.
Similar to the above discussions on the single region problems,
the error in the partial-HBIE (ii) is less
than that in the full-HBIE (iii) in the case of multi-regions problems.
Based on the comparison between model-(A) and either ({A:\scalebox{1.0}{$\,|\,$}}) or ({A:\scalebox{1.2}{$\times$}}),
we can observe that
the errors in the multi-media problems
with the interface boundaries (models-({A:\scalebox{1.0}{$\,|\,$}}) or ({A:\scalebox{1.2}{$\times$}})) are larger
than that in the single region problem (model-(A)).
One of the reasons could be the same reason for which the model, including the boundary
with $\overline{q}=0$ such as model-(B),
has a larger error than that with $\overline{u}=0$ such as model-(A), and
this is discussed in the description of \Fig{error_dist-open_end-1reg}.
The internal field is evaluated by \EQ{internal-field}
as the boundary integral where the boundary encloses the domain considered.
In the case of the multi-regions partitioned by the interface boundaries,
the boundary enclosing a single region must include the continuous boundary where $u\ne0$.
Therefore,
the error in $h_{i'}^j$ contributing from the interface boundary in \EQ{internal-field}
is added to the total error in the multi-region problems;
whereas it is not added from the boundary with $\overline{u}=0$ in model-(A).
The other reason is a difference in the distances between the internal points and their nearest boundary;
an average of distances in the multi-region problem is shorter than that of the single region problem.
Since the contributions $h_{i'}^j$ and $g_{i'}^j$ increase with decreasing distance,
the error in the multi-region problems is larger than that in the single-region problem.
The difference in error between the multi- and single-region problem can also be found in the comparison between model-(B) and ({B:\scalebox{1.2}{$\times$}}).

Consequently,
we can demonstrate that the formulation based on the HBIE is applicable
without rank deficiency even in the cases involving the interface boundaries,
which is similar to the corner nodes in the single region problems.

One may question whether the partial-HBIE method can avoid spurious solutions
of an external problem shown in Sec.~\ref{sec:intro}
since the partial-HBIE method uses both CBIEs and HBIEs
like the Burton-Miller method does to avoid spurious solutions.
In the Burton-Miller method, a linear combination of CBIEs and HBIEs with an appropriate combination factor
is used to prevent spurious solutions from arising when some of the sub-matrices become singular.
Whereas, in the partial-HBIE method, the matrix equation consists of two sets of equations, CBIEs and HBIEs, without any modifications.
Each set of equations has sub-matrices which may potentially produce spurious solutions.
The partial-HBIE method, therefore, cannot avoid spurious solutions.
To avoid them, we should use other methods,
such as the Burton-Miller method \cite{Burton-Miller:1971}, CHIEF \cite{Schenck:1968},
or a virtual boundary method \cite{Tomioka:1993,Tomioka:1994} which divides the external region into multiple regions by virtual boundaries to stop the external region from surrounding the internal region.
If we use a virtual boundary method, the issue of multiple-duplicated nodes has to be solved;
that is, however, not very difficult when using the method proposed in this paper.

\section{Conclusion}
\label{sec:conclusion}
The method of using double nodes at corners is a useful approach to uniquely
define the normal direction.
However, a set of simultaneous equations in CBIE formulation produces rank deficient problems
in the following cases:
both sub-nodes belonging to any double node are imposed by Dirichlet conditions;
an intersection of the interface boundary located between different media
is not connected to the boundary imposed by ordinary boundary conditions;
and an interface boundary is connected to the two boundaries imposed by Dirichlet conditions.
This means that the applicable problem that uses the double nodes are limited in the CBIE formulation.

In contrast,
when the coefficient matrix is constructed by HBIEs,
the rank is not reduced
for any combination of boundary conditions, including interface conditions.
However, the contribution coefficients between nodes in HBIEs are less accurate
than those in CBIEs for the problem without rank deficiency because of two reasons;
a HBIE exhibits a stronger singularity of the integrand than a CBIE,
and most of the coefficients are multiplied by the dyadic tensor with a large norm.
Furthermore,
the computational cost of evaluating the coefficients of HBIEs is higher than that of CBIEs.

To address the rank deficiency problem in CBIEs and the drawbacks in HBIEs,
the coupling approach presented in this paper called partial-HBIE is the best choice.
In the partial-HBIE, most node equations are constructed by CBIEs,
and only the sub-node equations related to corners are constructed by HBIEs.

The method that uses HBIEs demonstrates the following advantages compared to other methods:
it does not require
any additional local relation between nodal points around double nodes,
any extra boundary integral equation,
and it does not require a least-square method, which can be computationally time-consuming.
Furthermore,
the partial-HBIE can be applied by only switching the sub-node equation
for the double node from a CBIE to a HBIE;
therefore, we can be relieved from the efforts involved in preparing input data and complex coding.

%%%%%%%%%%%%%%%%%%%%%%%%%%%%%%%%%%%%%%%%%%%%%%%%%%%%%%%%%%%%%%%%%%

\section*{Acknowledgments}

This work was supported by JSPS KAKENHI Grant Number 18K04158.

%%%%%%%%%%%%%%%%%%%%%%%%%%%%%%%%%%%%%%%%%%%%%%%%%%%%%%%%%%%%%%%%%%%%

\bibliographystyle{elsarticle-num}
\bibliography{paper}

%%%%%%%%%%%%%%%%%%%%%%%%%%%%%%%%%%%%%%%%%%%%%

\newpage

\end{document}